\documentclass[twoside, 12pt]{article}
\usepackage[latin1]{inputenc}
\usepackage{amscd}
\usepackage{amsmath}
\usepackage{amsfonts}
\usepackage{amssymb}
\usepackage{amsthm}
\usepackage{comment}
\usepackage[francais,english]{babel}
\usepackage[matrix,arrow]{xy}

\renewcommand{\tilde}{\widetilde}
\renewcommand{\hat}{\widehat}

\newcommand{\norm}[1]{\left\| #1 \right\|}
    \DeclareMathOperator{\Ree}{Re}

    \DeclareMathOperator{\sgn}{sgn}
    
    \DeclareMathOperator{\Hom}{Hom}

    \DeclareMathOperator{\Cor}{Cor}
    \DeclareMathOperator{\Leb}{Leb}
    \DeclareMathOperator{\dLeb}{dLeb}

\newcommand{\de}{{\rm d}}
\newcommand{\dd}{\, {\rm d}}

\renewcommand{\geq}{\geqslant}
\renewcommand{\leq}{\leqslant}
\newcommand{\T}{\mathbb{T}}

\newcommand{\N}{\mathbb{N}}
\newcommand{\Z}{\mathbb{Z}}
\newcommand{\R}{\mathbb{R}}

\newcommand{\C}{\mathbb{C}}

\renewcommand{\phi}{\varphi}
\renewcommand{\epsilon}{\varepsilon}
\newcommand{\tq}{\ |\ }
\newcommand{\D}{\mathbb{D}}
\newcommand{\boL}{\mathcal{L}}

\newcommand{\boB}{\mathcal{B}}
\newcommand{\boN}{\mathcal{N}}
\newcommand{\boBB}{\mathcal{A}}

\newtheorem{thm}{Theorem}[section]
\newtheorem*{thm*}{Theorem}
\newtheorem{prop}[thm]{Proposition}
\newtheorem{lem}[thm]{Lemma}
\newtheorem{cor}[thm]{Corollary}

\theoremstyle{definition}

\newtheorem*{rmq}{Remark}

\title{Central limit theorem and stable laws for intermittent maps
  \footnote{\emph{keywords}: decay of correlations, intermittency, countable
    Markov shift, central limit theorem, stable laws, Wiener's
  Lemma.
  \emph{2000 Mathematics Subject Classification:} 37A30, 37A50, 37C30,
  37E05, 47A56, 60F05.}}
\author{Sébastien Gouëzel
  \footnote{Département de Mathématiques et Applications, 
École Normale Supérieure, 45 rue d'Ulm 75005 Paris
(France). e-mail \texttt{Sebastien.Gouezel@ens.fr}}}
\date{November 2002}
\begin{document}

\maketitle

\begin{abstract}
In the setting of abstract Markov maps, we prove results concerning the
convergence of renormalized Birkhoff sums to normal laws or stable
laws. They apply to one-dimensional maps with a neutral fixed
point at $0$ of the form $x+x^{1+\alpha}$, for $\alpha\in (0,1)$. In
particular, for $\alpha>1/2$, we show that the Birkhoff sums of a Hölder
observable $f$ converge to a normal law or a stable law, depending on
whether $f(0)=0$ or $f(0)\not=0$. The proof uses spectral techniques
introduced by Sarig, and Wiener's Lemma in non-commutative Banach algebras.
\end{abstract}

\section{Introduction and statement of results}

\subsection{Introduction}
Recently, general methods have been devised to prove the Central
Limit Theorem for the Birkhoff averages of mixing dynamical
systems. More precisely, if $T:X\to X$ is a map from a
space to itself preserving a probability measure $\nu$, and $f:X\to \R$ is a
function with zero average, we say that the Central Limit Theorem
(CLT) holds for $(T,\nu,f)$ if
  \begin{equation*}
  \frac{S_nf}{\sqrt{n}}=
  \frac{1}{\sqrt{n}} \sum_{i=0}^{n-1} f\circ T^i \to \boN(0,\sigma²),
  \end{equation*}
the convergence being in distribution (with respect to the measure
$\nu$).

To prove this kind of results, the first step is often to show
that the correlations $\Cor(f,f\circ T^n)=\int f\cdot f\circ T^n
\dd\nu$ are summable. Then one may use this to construct a reverse
martingale, for which the CLT is known to be true. Clearly, when
the correlations are not summable, this method does not work.

In this paper, we will be interested in cases where $T$ has a
neutral fixed point (or in the Markovian analogues of this
situation), where the classical methods to estimate the decay of
correlations do not work (however, Hu recently announced a proof of 
the
summability of the correlations \cite{Hu:annonce}). We will prove
that in some cases, the CLT still holds. Moreover, in situations
where it does not hold, we will show that there are other limit
laws appearing, namely stable laws, even when the function $f$ is
arbitrarily smooth, for example $C^\infty$ or analytic. This kind
of result has been announced by Fisher and Lopes 
(see \cite{fisher_lopes:noise} and references therein) for piecewise linear 
maps.

Our main results (Theorems \ref{TCL_markov} and
\ref{thm_loi_stable_markov}) will be stated in the context of abstract
Markov maps.
Applications to one-dimensional maps are then given:
if a map has a neutral fixed point at $0$ of the form $x+x^{1+\alpha}$
with $\alpha\in (1/2,1)$ and is expanding elsewhere,
then the Birkhoff sums of a Hölder observable $f$
converge to a normal law or a stable law depending on whether $f(0)=0$
or $f(0)\not=0$. Applications are also given for other types of
neutral fixed points, and for unbounded observables, for example
$x^{-\beta}$, leading to new results even in the case $\alpha\in
(0,1/2)$ (see also \cite{raugi:dim1} in this case, where he studies
the decay of correlations for non Hölder observables).

The philosophy of the method is that, if it is possible to induce on a subset
$Y$ of the space so that the induced map is uniformly expanding,
then the behavior of the Birkhoff sums $\sum_0^{n-1}f(T^i x)$ is
in fact dictated by the behavior on $Y$ of the Birkhoff sums
$\sum_0^{p-1}f_Y(T_Y^i x)$, where $T_Y=T^\phi$ is the induced map on
$Y$ ($\phi$ is the first return time), and
$f_Y(x)=\sum_0^{\phi(x)-1} f(T^i x)$. As convergence to normal laws
or stable laws is known in the uniformly expanding case (see
\cite{aaronson_denker} for example), we get the same kind of
results in our nonuniform setting.

\subsection{The result for Markov maps}
\subsubsection{Definition of Markov maps}
\label{section_def_markov_maps}
For these definitions, see for example
\cite{aaronson_denker_urbanski}, 
\cite{aaronson:book} or \cite{sarig:decay}.

A Markov map is a non-singular transformation $T$ of a Lebesgue
space $(X,\boB,m)$ of finite mass, together with a measurable partition $\alpha$
of $X$ such that if $a\in \alpha$ then $m(a)>0$, $Ta$ is a union (mod
$m$) of elements of $\alpha$, and $T:a\to Ta$ is invertible.
Moreover, it is assumed that the completion of $\bigvee_0^\infty
T^{-i}\alpha$ is $\boB$, i.e. the partition separates the points.

A Markov map is said to be \emph{topologically mixing} if
$\forall a,b \in \alpha, \exists N, \forall n \geq N, b\subset
T^n a$. This corresponds to topological mixing for the
topology defined by the cylinders
$[a_0,\ldots,a_{n-1}]=\bigcap_{i=0}^{n-1}T^{-i}a_i$ (where
$a_0,\ldots,a_{n-1}\in \alpha$).

If $\gamma \subset \alpha$ is a subset of the partition $\alpha$,
and $Y=\bigcup_{a\in \gamma} a$, the induced map $T_Y:Y\to Y$ is
defined as the first return map from $Y$ to $Y$, i.e.
$T_Y=T^{\phi_Y}$, where $\phi_Y(x)=\inf\{ n\geq 1\tq T^n(x)\in
Y\}$ is the return time to $Y$.
Let $\delta=\{[a,\xi_1,\ldots,\xi_{n-1},Y] \tq
a\in \gamma, \xi_1,\ldots,\xi_{n-1}\in \alpha-\gamma,
[a,\xi_1,\ldots,\xi_{n-1},Y]\not=\emptyset\}$: this is a
partition of $Y$, for which $T_Y$ is a Markov map (with the
measure $m_Y=m_{|Y}$). The cylinders for this partition will be
denoted by $[d_0,\ldots,d_{l-1}]_Y$ (with $d_0,\ldots,d_{l-1}\in
\delta$). If $d=[a,\xi_1,\ldots,\xi_{n-1},Y]\in \delta$, its
image is $T_Y d=T \xi_{n-1}\cap Y$ -- hence, it is
$\gamma$-measurable.

For $x,y\in Y$, we will write $s(x,y)$ for the separation time of
$x$ and $y$ under $T_Y$, i.e. $s(x,y)=\inf \{n\in \N \tq T_Y^n
x\text{ and }T_Y^n y \text{ are not in the same element of
}\delta\}$. This separation time is extended to the whole space
$X$: for $x,y\in X$, let $x'$ and $y'$ be their first returns to
$Y$. If $T^ix$ and $T^iy$ stay in the same element of the
partition $\alpha$ until the first return to $Y$, set
$s(x,y)=s(x',y')+1$. Otherwise, $s(x,y)=0$.

If $\theta<1$, we will say that a function $f$ on $Y$ (resp. $X$) is locally
$\theta$-Hölder if there exists a constant $C$ such that
$\forall x,y\in Y$ (resp. $X$) with $s(x,y)\geq 1$, we have $|f(y)-f(x)|\leq C
\theta^{s(x,y)}$. The smallest such constant is called the
$\theta$-Hölder constant of $f$. In fact, slightly abusing notation,
we shall say that a
measurable function is locally $\theta$-H\"older if there exists 
a 
locally $\theta$-H\"older version of this function.

The transfer operator $\hat{T_Y}$ associated to $T_Y$ and
acting on integrable functions is defined
by $\int f\cdot (g\circ T_Y) \dd m_Y=\int (\hat{T_Y}f) \cdot g \dd
m_Y$. It can be written as
$\hat{T_Y}f(x)=\sum_{T_Y y=x}g_{m_Y}(y)f(y)$, where the weight
$g_{m_Y}$ is defined by $g_{m_Y}=\frac{\de m_Y}{\de m_Y\circ
T_Y}$, i.e. it is the  inverse of the jacobian of $T_Y$.
We will say that the distortion is locally Hölder if
$\log g_{m_Y}$
is locally $\theta$-Hölder for some $\theta<1$.

Finally, we will say that $T_Y$ has the \emph{big image property}
if there exists a constant $\eta>0$ such that $\forall d\in
\delta, m[T_Y d]\geq \eta$. This is in particular the case when
the partition $\gamma$ is finite.

\subsubsection{Main results}

For the following theorems, $(X,\boB,T,m,\alpha)$ will be a
topologically mixing probability preserving Markov map. We assume
that $\gamma \subset \alpha$ is such that the induced map $T_Y$ on
$Y=\bigcup_{a\in \gamma}a$ has a locally $\theta$-Hölder distortion
for some $\theta$,
and that it has the big image property. The function $\phi=\phi_Y$ will denote
the return time from $Y$ to itself.

We will write $I(X,Y)$ for the set of functions $f:X\to \R$ which induce
``nice'' functions on $Y$. More precisely, $f\in I(X,Y)$ if
  \begin{enumerate}
\item
  $\forall N\in \N$, there exist constants $C_N$ and $0<\theta_N<1$ such
  that, on $Y_N=\bigcup_0^{N-1} T^i(Y)$, we have $|f|\leq C_N$ and, if
  $x,y\in Y_N$ satisfy $s(x,y)\geq 1$, then $|f(x)-f(y)|\leq C_N
  \theta_N^{s(x,y)}$.
\item
  Writing $f_Y(x)=\sum_{i=0}^{\phi(x)-1} f(T^i x)$ for the map induced by $f$
  on $Y$, there exists $\theta<1$ such that
  $\sum_{d\in\delta} \phi(d) m[d] D_\theta
  f_Y(d)<\infty$, where $D_\theta f_Y(d)$ is the least $\theta$-Hölder
  constant of $f_Y$ restricted to the element $d$ of the partition
  $\delta$.
  \end{enumerate}
Note that, increasing $\theta$ if necessary, we can assume that the
distortion of $T_Y$ is $\theta$-H\"older and that 
$\sum_{d\in\delta} \phi(d) m[d] D_\theta f_Y(d)<\infty$, for the same
$\theta<1$.

For example, if $f$ is bounded and locally $\theta$-Hölder, and
the induced map $f_Y$ is also locally $\theta$-Hölder, then these
conditions are satisfied (note that $\sum \phi(d)m[d]=m(X)<\infty$, by
Kac's Formula). The weaker conditions above are useful in
the applications, for example to unbounded
observables.

\begin{thm}[Central Limit Theorem]
\label{TCL_markov}
Let $f\in I(X,Y)$ be integrable with $\int f\dd m=0$ and let the
induced function $f_Y$ be in $L^2(\dd m_Y)$.

Then there exists a constant $\sigma^2\geq 0$ such that
$\frac{1}{\sqrt{n}}\sum_{i=0}^{n-1} f\circ T^i$ converges in
distribution to a normal law $\boN(0,\sigma^2)$. Moreover,
$\sigma^2=0$ if and only if there exists a measurable function
$\chi:X\to \R$ satisfying $f=\chi \circ T-\chi$.
\end{thm}

For the next theorem, we will need the notion of stable law. The
law of the random variable $X$ is said to be stable if there
exist i.i.d. random variables $X_k$ and constants $A_n\in \R$ and
$B_n>0$ such that
  \begin{equation*}
  \frac{\sum_{i=0}^{n-1} X_i -A_n}{B_n}\to X \text{ in distribution.}
  \end{equation*}
In this case, we say that $X_1$ is attracted to $X$.
The stable laws are completely classified (see for example
\cite{feller:2}) and depend on many parameters, the most important
of which is certainly the \emph{index} $p\in (0,2]$. The case $p=2$
corresponds to the normal law, while if $p\in (1,2)$ the
stable laws are in $L^1$ but not in $L^2$. In fact, for $p\in
(0,1)\cup (1,2)$,
the index can be characterized as follows: there exist constants
$c_1\geq 0$ and $c_2\geq 0$ with $c_1 +c_2>0$, satisfying $P[X>x]
=(c_1+o(1)) x^{-p}$ and $P[X<-x]=(c_2+o(1))x^{-p}$. The case $p=1$ is
problematic and will not be included in the following discussion.

Following \cite{aaronson_denker}, for $p\in (0,1)\cup (1,2]$, $c>0$ and
$\beta \in [-1,1]$, we will denote by $X_{p,c,\beta}$ the law whose
characteristic function is
  \begin{equation*}
  E(e^{it X_{p,c,\beta}})=
  e^{-c|t|^p
  \left(1-i\beta \sgn(t) \tan\left(\frac{p\pi}{2}\right)\right)} .
  \end{equation*}
This is indeed a stable law of index $p$. When $p\in (1,2]$, it is centered,
i.e. $X_{p,c,\beta}$ is of zero expectation. Note that, when $p=2$,
the value of $\beta$ is irrelevant.

We recall also that a function
$L:(0,\infty)\to \R$ is slowly varying if $\forall x>0,
\lim_{y\to\infty}
 \frac{L(yx)}{L(y)}=1$.
Then the random variables $Y$ attracted to stable
laws of index $p$ can be described as follows (\cite{feller:2}):
\begin{itemize}
\item
For $p\in (0,1)\cup (1,2)$, $P[Y>x]=(c_1+o(1))x^{-p}L(x)$ and
$P[Y<-x]=(c_2+o(1))x^{-p}L(x)$ for some $c_1,c_2\geq 0$ with
$c_1+c_2>0$, and a slowly varying function $L$.
\item
For $p=2$, either $Y\in L^2$, or $P[Y>x]=(c_1+o(1))x^{-2}l(x)$ and
$P[Y<-x]=(c_2+o(1))x^{-2} l(x)$, for constants $c_1,c_2\geq 0$ with
$c_1+c_2>0$, and a slowly varying function $l$ such that $u^{-1}l(u)$
is not integrable at $+\infty$.
In this case, we will write $L(x)=\int_{-x}^x u^2
\dd P_Y(u)$, which is an unbounded slowly varying function.
\end{itemize}
When $Y\not\in L^2$,
if the normalizing constants $B_n$ and $A_n$ are defined by
  \begin{equation*}
  nL(B_n)=B_n^p,\ \  A_n=\left\{ \begin{array}{ll}
  0& 0<p<1\\
  nE(Y)& 1<p\leq 2\\
  \end{array}\right. ,
  \end{equation*}
then, writing $Y_0,Y_1,\ldots$ for independent random variables having the
same distribution as $Y$, we have the convergence
  \begin{equation*}
  \frac{\sum_{i=0}^{n-1}Y_i -A_n}{B_n} \to X_{p,c,\beta}
  \end{equation*}
where
  \begin{equation*}
  c=\left\{ \begin{array}{ll}
  (c_1+c_2)\Gamma(1-p)\cos\left(\frac{p\pi}{2}\right) & p\in (0,1)\cup (1,2) \\
  \frac{1}{2} & p=2
  \end{array}\right.,\ \
  \beta=\frac{c_1-c_2}{c_1+c_2}.
  \end{equation*}

In the deterministic case, if $f\in I(X,Y)$, we can consider the
function $f_Y$ on $Y$. Since $m(Y)<1$, its distribution $\mu$ (defined
by $\mu(A)=m\{ y\in Y \tq f_Y(y)\in A\}$) is not a probability measure
on $\R$, but we can replace it by the modified distribution
$\tilde{\mu}=\mu+(1-m[Y])\delta_0$. We will say with a little abuse of
notation that the distribution of $f_Y$ is attracted to a stable law if
the distribution $\tilde{\mu}$ is in the domain of attraction of this
stable law.

\begin{thm}[Convergence to stable laws]
\label{thm_loi_stable_markov}
Let $f\in I(X,Y)$. Assume that
the induced function $f_Y$ on $Y$ is not in $L^2$, and that its
distribution is in the domain of attraction of a stable law of index
$p\in (0,1)\cup (1,2]$, as described previously.

Then $\frac{\sum_{i=0}^{n-1} f\circ T^i-A_n}{B_n}$ converges in
distribution to the stable law $X_{p,c,\beta}$ where
$A_n,B_n,p,c$ and $\beta$ are as in the random case above.
\end{thm}

Note that, in Theorem
\ref{TCL_markov}, the constant $\sigma^2$ depends on the interactions
at different times between $f$ and $f\circ T^n$, while in Theorem
\ref{thm_loi_stable_markov}, the limit distribution depends only on the initial
distribution (of $f_Y$).

\subsection{Application to one-dimensional maps}

Our main theorems can be applied to one-dimensional maps with a
neutral fixed point with prescribed behavior. For the sake of
simplicity, we will restrict ourselves to a specific case, but the
results can in fact be proved for a wide class of maps.

For $0<\alpha<1$, we consider the map from $[0,1]$ to itself introduced in
\cite{liverani_saussol_vaienti} (this paper contains also
historical references) and defined by
  \begin{equation}
  \label{definition_LSV}
   T(x)=\left\{ \begin{array}{cl}
   x(1+2^\alpha x^\alpha) &\text{if }0\leq x\leq 1/2
   \\
   2x-1 &\text{if }1/2<x\leq 1.
   \end{array}\right.
   \end{equation}
This map has a neutral fixed point at $0$, and is expanding
elsewhere. It is known that it has a unique absolutely continuous
invariant measure $\dd m$, whose density $h$ is Lipschitz on any interval
of the form $(\epsilon,1]$ (\cite[Lemma
2.3]{liverani_saussol_vaienti}).

For Hölder observables on $[0,1]$, 
Young (\cite{lsyoung:recurrence}) showed that
$\Cor(f,g\circ T^n)=\int f\cdot g\circ T^n \dd m-\int fg \dd m=
O(1/n^{1/\alpha-1})$. This rate is in fact optimal for a wide class of
functions (\cite{sarig:decay}). It implies a central limit theorem
when $\alpha<1/2$.

For $\alpha>1/2$, the above rate is not summable. However, extending
the methods of Sarig, we showed in \cite{gouezel:decay}
that, for functions $f$ vanishing in a neighborhood of the fixed point and
with zero integral, the correlations decay like $O(n^{-1/\alpha})$,
which gives again a central limit theorem. For functions with zero
average and equal to a nonzero constant in a neighborhood of the
fixed point, the correlations decay
exactly like $C/n^{1/\alpha-1}$ (\cite{gouezel:decay}), which is not
summable.

Our main theorems yield (see also a list of other consequences after
the proof of Theorem \ref{TCL_LSV}):

\begin{thm}
\label{TCL_LSV} Let $1/2<\alpha<1$, and $T:[0,1]\to [0,1]$ be the
corresponding Liverani-Saussol-Vaienti map \eqref{definition_LSV}. 
Let $h$ be the density of
its unique absolutely continuous invariant probability $\dd m$.
Let $f:[0,1]\to \R$ be Hölder, and
$\gamma$-Hölder on a neighborhood of $0$ with
$\gamma>\frac{2\alpha-1}{\alpha+1}$. Assume also $\int f\dd m=0$.

Then, if $f(0)=0$, there exists a constant $\sigma^2\geq 0$ such that
$\frac{S_n f}{\sqrt{n}}$ tends in distribution to
$\boN(0,\sigma^2)$. Moreover, $\sigma^2=0$ if and only if there exists
a measurable function $\chi$ with $f=\chi\circ T-\chi$.

If $f(0)\not=0$, then $\frac{S_nf}{n^\alpha}$ converges in distribution to the
stable law $X_{1/\alpha,c,\sgn(f(0))}$ with $c=\frac{h(1/2)}{4(\alpha/|f(0)|)^{1/\alpha}}
\Gamma(1-1/\alpha)\cos\left(\frac{\pi}{2\alpha}\right)$.
\end{thm}

This result does not address the question of the summability of
correlations for a function $f$ with $f(0)=0$. Hu announced recently
this summability (\cite{Hu:annonce}), which gives another proof of the central limit theorem.

\begin{proof}
We will assume that $f$ is $\gamma$-Hölder on the whole interval
$[0,1]$, to avoid notational problems.

We construct a Markov structure, starting from $x_0=1$ and setting
$x_{n+1}=T^{-1}(x_n)\cap [0,1/2]$. The map $T$ is then Markov for
the partition composed of the $(x_{n+1},x_n), n\in \N$. We will
induce on $Y=(1/2,1)$. The distortion control on the induced map
is classical (see for example \cite{liverani_saussol_vaienti}), it
remains to study the properties of $f$ and $f_Y$. We will write
$s(x,y)$ for the separation time of $x,y\in [0,1]$ under the
induced map $T_Y$.

As the induced map on $Y$ is uniformly expanding, with a factor
$\lambda>1$, we obtain $|x-y|\leq C \lambda^{-s(x,y)}$. As $f$ is
$\gamma$-Hölder, we get
$|f(x)-f(y)|\leq D |x-y|^\gamma \leq E
\left(\lambda^{-\gamma}\right)^{s(x,y)}$, with
$\theta=\lambda^{-\gamma}<1$. Thus, $f$ is bounded and locally Hölder.

To prove that $f_Y$ is locally Hölder for the separation time $s$,
the argument is essentially
the same, with some distortion. The induced partition for which the
induced map $T_Y$ is Markov is composed of the intervals $J_k$ that
are mapped by $T$ onto $(x_k,x_{k-1})$. We know that
the distortion of $T^k:(x_{k+1},x_k)\to (1/2,1)$ is bounded
(see \cite{liverani_saussol_vaienti}), that $x_k\sim
C/k^{1/\alpha}$ for $C=\frac{1}{2}\alpha^{-1/\alpha}$, and that
$|x_k-x_{k+1}|=2^\alpha x_{k+1}^{\alpha+1} \sim K/k^{1+1/\alpha}$. Thus,
on $(x_{k+1},x_k)$, we have
$(T^k)'\asymp |T^k(x_{k+1},x_k)|/|(x_{k+1},x_k)|\asymp k^{1+1/\alpha}$
(the symbol $a \asymp b$ means here that the ratio $a/b$ is bounded
from above and from below). We can assume without loss of generality
that $\gamma(1+1/\alpha)<1$.
Then, for $x,y\in J_k$,
  \begin{align*}
  |f_Y(x)-f_Y(y)|&
  \leq|f(x)-f(y)|+\sum_{i=1}^{k-1}C |T^i x-T^i
  y|^\gamma
  \\&
  \asymp |f(x)-f(y)|+\sum_{i=1}^{k-1} \left(\frac{|T^kx-T^k
    y|}{(T^{k-i})'_{|(x_{k-i+1},x_{k-i})}}\right)^\gamma
  \\&
  \asymp |f(x)-f(y)|+|T^k x-T^k y|^\gamma
  \sum_{i=1}^{k-1} \left(\frac{1}{(k-i)^{1+1/\alpha}}\right)^\gamma
  \\&
  \leq K \theta^{s(x,y)} k^{1-\gamma(1+1/\alpha)}
    \end{align*}
Thus,
  \begin{equation*}
  \sum \phi(J_k) m[J_k] D_\theta f_Y(J_k)
  \leq K \sum k \frac{1}{k^{1+1/\alpha}} k^{1-\gamma(1+1/\alpha)}
  \end{equation*}
which is finite when $\gamma >\frac{2\alpha-1}{\alpha+1}$. This proves
that $f\in I(X,Y)$.

We assume now that $f(0)=0$ and show that $f_Y\in L^2(Y,\dd m)$. We may
suppose that $\gamma<\alpha$. As $\dd m$
and $\dLeb$ are equivalent on $Y$, it is sufficient to show that $f_Y\in
L^2(Y,\dLeb)$. If $x\in J_n$, then for $1\leq i\leq {n-1}$, $T^i x \in
(x_{n+1-i},x_{n-i})$, whence $|f(T^i x)|\leq C |T^i x|^\gamma\leq C
|x_{n-i}|^\gamma$. This implies that $|f_Y(x)|\leq |f(x)|+C
\sum_{j=1}^{n-1} |x_j|^\gamma \leq C+C\sum_{j=1}^{n-1}
\frac{1}{j^{\gamma/\alpha}}$.
Hence, $|f_Y|\leq D n^{1-\gamma/\alpha}$ on $J_n$, with
$\Leb(J_n) \asymp x_n-x_{n+1} \asymp 1/n^{1+1/\alpha}$. Thus, $\int
|f_Y|^2 \leq C \sum \frac{n^{2-2\gamma/\alpha}}{n^{1+1/\alpha}}$. This
is summable when $\gamma>\alpha-\frac{1}{2}$, which is implied by
$\gamma>\frac{2\alpha-1}{\alpha+1}$. We have proved that $f_Y\in L^2$,
whence we can use Theorem \ref{TCL_markov} to conclude.

We assume finally that $f(0)>0$ (the case $f(0)<0$ is analogous).
Writing $f=f(0)+\tilde{f}$, the previous
estimate shows that $\tilde{f_Y}(x)
=\sum_0^{\phi(x)-1} \tilde{f}(T^i x)\in L^2(Y)$. Writing $u(x)=kf(0)$ on
$J_k$, we have $f_Y=\tilde{f_Y}+u$ on $Y$, and
  \begin{equation*}
  m[u>kf(0)]
  =m\left[\bigcup_{k+1}^\infty J_l\right]=m(1/2,x_k/2+1/2)
  \sim h(1/2)\frac{x_k}{2}
  \sim h(1/2)\frac{1}{4}(\alpha k)^{-1/\alpha}.
  \end{equation*}
Thus, we obtain $m[u>x]\sim \frac{h(1/2)}{4 (\alpha
  x/f(0))^{1/\alpha}}=\frac{c_1}{x^{1/\alpha}}$.
Adding $\tilde{f_Y}\in L^2$ does not change the asymptotics since
  $m[|\tilde{f_Y}|>x]=o(x^{-2})$ because $\tilde{f_Y}\in L^2$.
We have proved that $m[f_Y>x]=(c_1+o(1))x^{-p}$ for $p=1/\alpha$. 
In the same way,
we get $m[f_Y<-x]=o(x^{-p})$. This enables us to use Theorem
\ref{thm_loi_stable_markov}, which concludes the proof.
\end{proof}

We enumerate other consequences of Theorems \ref{TCL_markov} and
\ref{thm_loi_stable_markov} in other one-dimensional
cases:
\begin{enumerate}
\item
For $\alpha\in (0,1/2)$, the function $f_Y$ is always in $L^2$ and the
condition $\gamma >\frac{2\alpha-1}{\alpha+1}$ becomes empty. Thus, we
get another proof of the classical central limit theorem, for any
Hölder function $f$.

\item
For $\alpha=1/2$, the condition on $\gamma$ is also empty, and we get
the following result for a Hölder $f$ with zero average:
\begin{itemize}
\item
If $f(0)=0$, there is a central limit theorem.
\item
If $f(0)\not=0$, for example $f(0)>0$, then we get as above
$m[f_Y>x]=(c_1+o(1))x^{-2}$ for $c_1=\frac{h(1/2) \sqrt{2 f(0)}}{4}$. We are in the
nonstandard domain of attraction of the normal law (i.e. $f_Y$ is
attracted to a normal law but is not in $L^2$), for $l(x)=1$,
whence $L(x)\sim \frac{c_1}{2}\log x$. The constant $B_n$ is then
$\frac{\sqrt{c_1}}{2}\sqrt{n \log
  n}$, and we get the convergence
  \begin{equation*}
  \frac{\sum_{i=0}^{n-1}f\circ T^i}{\frac{\sqrt{c_1}}{2}\sqrt{n\log n}} \to
  \boN(0,1).
  \end{equation*}
\end{itemize}

\item
For $\alpha\in (0,1)$, we can also look at $f(x)=x^{-\beta}+c$ with
$\beta>0$, where
$c$ is taken so that $\int f=0$. Note that $f$ is integrable with
respect to the invariant measure $m$ if and only if $\alpha+\beta<1$
(since the density of $m$ behaves like $x^{-\alpha}$ close to $0$). A
small computation shows that $D_\theta f_Y(J_n)\asymp n^{\beta/\alpha}$,
whence $\sum n m[J_n]D_\theta f_Y(J_n)<\infty$ as soon as
$\alpha+\beta<1$.
Finally, on $J_n$, $f_Y\sim C n^{\beta/\alpha+1}$, whence $m[f_Y>x]\sim C
x^{-\frac{1}{\alpha+\beta}}$. Thus, when $\alpha+\beta<1/2$, we obtain
a usual central limit theorem, while when $\alpha+\beta=1/2$ we have a
central limit theorem with normalization $\frac{1}{\sqrt{n\log n}}$,
and for $1/2<\alpha+\beta<1$ we obtain a convergence to a stable law
for the normalization $\frac{1}{n^{\alpha+\beta}}$. So, even when
$\alpha<1/2$, we can have convergence to a stable law if the
observable is unbounded.

Note that Theorems \ref{TCL_markov} and
\ref{thm_loi_stable_markov} can not be used directly in this case,
since the function $f$ is not bounded on $\bigcup_0^{N-1}
T^i(Y)$. However, what is important in the proof is that we can take
an increasing sequence $Z_q$ with $\bigcup Z_q=X$ such that the
induction on $Z_q$ has good properties. In this case, we can take
$Z_q=(x_q,1)$ and the proof is exactly the same. In an equivalent way,
we can construct an inverse tower $\bigcup_{0}^{N-1} T^{-i}(Y)$ and
work with this inverse tower, the proofs being almost identical.

\item
For fixed points such as $x+x^\alpha \log x$ with $\alpha\in (0,1)$,
we also get other slowly
varying functions, whence convergence with other normalizing
factors.
\end{enumerate}

\subsection{Strategy of the proof}
Results of convergence to normal laws or stable laws are already
known when it is possible to obtain a spectral gap for the
transfer operator, for example for uniformly expanding dynamics (see
for example \cite{guivarch-hardy} and \cite{aaronson_denker}).
In the following, this case will be referred to as the
\emph{uniform} case. Such results apply in particular to the
induced maps, and they are obtained by perturbing the transfer
operators.

The strategy is then to translate the results on the induced
map to the whole space. For this, we use \emph{first
return transfer operators}, which have been used by Sarig in
\cite{sarig:decay} (see also \cite{isola:ergodic_degree}).
We will perturb these operators, and transfer
the results in the uniform case to these first return transfer
operators, through the abstract Theorem \ref{TCLabstrait}
which is in fact the main technical part of the proof.

This gives information on the behavior of the Birkhoff sums, but
only on the part of the space on which we have induced. We then
have to induce on larger and larger parts of the space, to recover
more and more information, which will complete the proof.

In Section \ref{section:abstract_result}, we will prove the
abstract result referred to above, Theorem \ref{TCLabstrait}.
In Section
\ref{section:uniform_case}, we will establish estimates in the
uniform case in our context,
which is slightly more general than the results in the
literature. Finally, we will prove
the main theorems in Section \ref{section:proof_thm}.

\section{An abstract perturbative theorem}
\label{section:abstract_result}

In the following, $\D$ will denote $\{z \in \C \tq |z|<1\}$,
and $\overline{\D}=\{z\in \C \tq |z|\leq 1\}$.

\subsection{The result}
The goal of this section is to prove the following theorem:

\begin{thm}
\label{TCLabstrait} Let $\boL$ be a Banach space and $R_n\in
\Hom(\boL,\boL)$ be operators on $\boL$ with $\norm{R_n}\leq r_n$
for a sequence $r_n$ such that $a_n=\sum_{k>n}r_k$ is
summable. Write $R(z)=\sum R_n z^n$ for $z\in \overline{\D}$.
Assume that $1$ is a simple isolated eigenvalue of $R(1)$
and that $I-R(z)$ is invertible for $z\in \overline{\D}-\{1\}$.
Let $P$ denote the spectral projection of $R(1)$ for the
eigenvalue $1$, and assume that $PR'(1)P=\mu P$ for some $\mu>0$.

Let $R_n(t)$ be an operator depending on $t\in [-\alpha,\alpha]$,
continuous at $t=0$, with $R_n(0)=R_n$ and $\norm{R_n(t)}\leq C
r_n$ $\forall t\in[-\alpha,\alpha]$, for some constant $C>0$.
For $z\in \overline{\D}$ and $t\in
[-\alpha,\alpha]$, write
  \begin{equation*}
  R(z,t)=\sum_{n=1}^\infty z^n R_n(t).
  \end{equation*}
This is a continuous perturbation of $R(z)$. For $t$ small and $z$
close to $1$, $R(z,t)$ is close to $R(1)$, whence it
admits an eigenvalue $\lambda(z,t)$ close to $1$. Assume that
$\lambda(1,t)=1-(c+o(1))M(|t|)$ for $c\in \C$ with
$\Re(c)>0$, and some continuous function $M:\R_+\to \R_+$ vanishing
only at $0$. Then
  \begin{enumerate}
  \item There exists $\epsilon_0>0$ such that $\forall
|t|<\epsilon_0$, $I-R(z,t)$ is invertible for all $z\in\D$. We can
write $(I-R(z,t))^{-1}=\sum T_{n,t} z^n$.
  \item
  Furthermore, there
exist functions $\epsilon(t)$ and $\delta(n)$ tending to $0$ when
$t\to 0$ and $n\to\infty$ such that $\forall |t|<\epsilon_0$,
$\forall n\in \N^*$, we have
  \begin{equation}
  \label{EquationTCLabstrait}
  \norm{T_{n,t}-\frac{1}{\mu}\left(1-\frac{c}{\mu}M(|t|)\right)^n
  P} \leq \epsilon(t)+\delta(n).
  \end{equation}
  \end{enumerate}
\end{thm}

For example, when $M(t)=t^p$, this gives the convergence
$T_{n,t/n^{1/p}}\to \frac{e^{-c|t|^p/\mu}}{\mu}P$.

\begin{rmq}
There is an analogue of this theorem when the perturbation takes
place on $[0,\alpha]$ (with exactly the same proof). Thus, for a
perturbation where $\lambda(1,t)=1-a|t|^p+ib \sgn(t)
|t|^p+o(|t|^p)$ for $t\in [-\alpha,\alpha]$,  we get once again
the same result, using twice this theorem, on $[-\alpha,0]$ and
$[0,\alpha]$.
\end{rmq}

\subsection{Wiener's Lemma in Banach algebras}

The classical Wiener's Lemma says that if $f:S^1\to \C$ has
summable Fourier coefficients and is everywhere nonzero, then the
Fourier coefficients of $1/f$ are also summable (see for example
\cite{kahane:sommable}). We will use an analogue of this result for
functions taking their values in non-commutative Banach algebras.

Let $\boBB$ be a Banach algebra with unit (its norm will be written
$\norm{\ }_\boBB$). We will denote by $A(\boBB)$ the set of
continuous functions $f:S^1 \to \boBB$ whose Fourier coefficients,
defined by
  \begin{equation*}
  c_n(f)=\frac{1}{2\pi} \int_0^{2\pi}
  e^{-in\theta}f(e^{i\theta})\dd\theta,
  \end{equation*}
are summable in norm. Then $f$ coincides with the sum
$\sum_{n\in\Z} c_n(f) z^n$ of its Fourier series (see
\cite[Proposition 2.6]{gouezel:decay}). In particular, if $f,g\in
A(\boBB)$, then $c_n(fg)=\sum_{k\in \Z} c_{n+k}(f) c_{n-k}(g)$,
whence $fg \in A(\boBB)$. Moreover, $A(\boBB)$ is a Banach algebra
for the norm $\norm{f}_{A(\boBB)}= \sum \norm{c_n(f)}_\boBB$.

For $f:S^1\to \boBB$, we will also write
$\norm{f}_\boBB=\sup\{\norm{f(z)}_\boBB \tq z\in S^1\}$. We have of
course $\norm{f}_\boBB \leq \norm{f}_{A(\boBB)}$, but the other
inequality is false.

The following result can be found in \cite{wiener_non_commutatif} (I
thank Omri Sarig for this reference).
\begin{thm}
\label{wiener_classique} Let $f\in A(\boBB)$ be such that for all
$z\in S^1$, $f(z)$ is invertible in the Banach algebra $\boBB$.
Then, setting $g(z)=f(z)^{-1}$, we also have $g\in A(\boBB)$.
\end{thm}

Let $\Omega$ be an open subset of $\C$ and $h:\Omega\to\C$ an
holomorphic function. If $x\in \boBB$ is such that its spectrum
$\sigma(x)$ is included in $\Omega$, it is possible to define
$h(x)\in \boBB$ using the Cauchy formula, by
$h(x):=\frac{1}{2i\pi}\int_\gamma \frac{h(u)}{uI-x}\dd u$, where
$\gamma$ is a path around $\sigma(x)$ in $\Omega$. This is
independent of $\gamma$ (see e.g.
\cite[VII.3]{dunford_schwartz:1}). For example, if $h(z)=z^n$ for
$n\in \Z$, then $h(x)=x^n$ (for the multiplication in the algebra
$\boBB$). For another example, take $\boBB=\Hom(\boL,\boL)$ where
$\boL$ is a Banach space, and $x\in \boBB$ for which $1$ is
isolated in $\sigma(x)$. Then, if $h$ is equal to $1$ in a small
neighborhood of $1$, and to $0$ outside of this neighborhood, then
$h(x)$ is the spectral projection associated to $1$.

\begin{thm}
  \label{wiener_general}
Let $\Omega$ be an open subset of $\C$ and $h:\Omega\to\C$ be
holomorphic. Let $f\in A(\boBB)$ be such that $\forall z\in S^1$,
$\sigma(f(z)) \subset \Omega$. If $g(z)=h(f(z))$, then we also
have $g\in A(\boBB)$.
\end{thm}
The previous Wiener's Lemma is a particular case of this result,
for $h(z)=1/z$.
\begin{proof}
As the spectrum is semi-continuous, there exists a path $\gamma$
in $\Omega$ around $\sigma(f(z))$ for every $z\in S^1$. Thus,
$g(z)=\frac{1}{2i\pi}\int_\gamma \frac{h(u)}{uI-f(z)}\dd u$. For
$u\in\gamma$ fixed, $M_u : z\mapsto \frac{h(u)}{uI-f(z)}$ is in
$A(\boBB)$ by Theorem \ref{wiener_classique}. As the inversion is
continuous in the Banach algebra $A(\boBB)$, $u\mapsto M_u$ is
continuous on the compact $\gamma$. Since $g=\int_\gamma M_u \dd
u$, this concludes the proof.
\end{proof}

Note that this result applies in particular to spectral
projections.

\subsection{Study of $(I-R(z,t))^{-1}$ at $t=0$}
In this section, we will describe the behavior of the
coefficients of $(I-R(z))^{-1}$, and in particular show that they
tend to $\frac{1}{\mu}P$. We will use techniques introduced by
Sarig in \cite{sarig:decay}. The main technical difference is
that in the proof of Theorem \ref{sarig_first_main} (an analogue of his
``first main lemma''), we will replace his estimates using
$C^{1+\alpha}$ regularity by applications of Wiener's Lemma.

The following result is interesting in itself, and its proof will be
instrumental in the proof of Theorem \ref{TCLabstrait}.

\begin{thm}
\label{sarig_first_main}
Let $R_n$ be operators on a Banach space
$\boL$ such that $a_n=\sum_{k>n}\norm{R_k}$ is summable.
Write $R(z)=\sum_{n=1}^\infty R_n z^n$ and assume that $1$ is a simple
isolated eigenvalue of $R(1)$, while $I-R(z)$ is invertible for
$z\in \overline{\D}-\{1\}$. Write $P$ for the spectral projection of
$R(1)$ for the eigenvalue $1$, and assume also that $PR'(1)P=\mu P$ for a
nonzero constant $\mu$.

Then the function $\left(\frac{I-R(z)}{1-z}\right)^{-1}$ can be
continuously extended to $\overline{\D}$, and its Fourier
coefficients on $S^1$ are summable. Writing
$(I-R(z))^{-1}=\sum_n T_n z^n$, we have $\sum_n \norm{T_n-T_{n-1}}
<\infty$.
\end{thm}
\begin{proof}
Write $S(z)=\frac{I-R(z)}{1-z}$.

\emph{Step 1: $S(z)^{-1}$ can be continuously extended to
$\overline{\D}$.}

The problem is the extension to $1$. For $z\in\overline{\D}$ close
to $1$, $R(z)$ has a unique eigenvalue $\lambda(z)$ close to $1$.
We write $P(z)$ for the corresponding spectral projection, and
$Q(z)=I-P(z)$. For $z$ close to $1$ but different from $1$,
  \begin{equation}
  \label{S-1_converge_en_1}
  S(z)^{-1}=\frac{1-z}{1-\lambda(z)}P(z)+ (1-z) (I-R(z)Q(z))^{-1}
  Q(z).
  \end{equation}
As $I-R(z)Q(z)$ is everywhere invertible on $\overline{\D}$, 
the second term has a
continuous extension to $1$. As $(I-R(1))P(1)=0$,
  \begin{equation}
  \label{quotient_tend_vers_mu}
  \frac{1-\lambda(z)}{1-z}P(z)
  =P(z)\frac{R(1)-R(z)}{1-z}P(z)+
  P(z)(I-R(1))\frac{P(z)-P(1)}{1-z}.
  \end{equation}
Note that
  \begin{equation*}
  \frac{R(z)-R(1)}{z-1}=\sum_{n=0}^\infty \left(\sum_{k=n+1}^\infty
  R_k\right)z^n.
  \end{equation*}
As $\sum_{k>n}\norm{R_k}$ is a summable sequence, the above sum converges
in norm, which implies
that its limit at $1$, denoted by $R'(1)$, is well
defined.

If $\delta>0$ is small enough, then
  \begin{equation*}
  P(z)=\frac{1}{2i\pi}\int_{|u-1|=\delta} \frac{1}{uI-R(z)}\dd u.
  \end{equation*}
Thus,
  \begin{equation}
  \label{pz-p1_sommable}
  \frac{P(z)-P(1)}{z-1}=\frac{1}{2i\pi}\int_{|u-1|=\delta}
  \frac{1}{uI-R(z)} \frac{R(z)-R(1)}{z-1} \frac{1}{uI-R(1)} \dd u.
  \end{equation}
As $\frac{R(z)-R(1)}{z-1}$ converges at $1$, we get the
convergence of $\frac{P(z)-P(1)}{z-1}$ at $1$, to a limit $P'(1)$.

In Equation \eqref{quotient_tend_vers_mu}, when $z$ tends to $1$,
the right hand term tends to
$P(1)R'(1)P(1)-P(1)(I-R(1))P'(1)=P(1)R'(1)P(1)$ since
$P(1)(I-R(1))=0$. Recall that $P(1)R'(1)P(1)=\mu P$ for some
nonzero $\mu$.

Let $\xi$ be a bounded linear functional on $\Hom(\boL,\boL)$ such that
$\xi(P(1))=1$. We apply $\xi$ to both members of Equation
\eqref{quotient_tend_vers_mu} and let $z$ go to $1$, which gives
that
  \begin{equation}
  \label{converge_vers_mu}
  \frac{1-\lambda(z)}{1-z}\to \mu.
  \end{equation}
As $\mu\not=0$, this
shows that, in Equation \eqref{S-1_converge_en_1}, both terms on
the right hand side converge at $1$.

\emph{Step 2: Construction of a function $\tilde{R}(z)$ on $S^1$,
  coinciding with $R(z)$ for $z$ close to $1$, such that it has
  everywhere a simple eigenvalue $\tilde{\lambda}(z)$ close to $1$,
  but different from $1$ if $z\not=1$,
  and with $\frac{\tilde{R}(z)-\tilde{R}(1)}{z-1}\in A(\boBB)$.}

We construct $\tilde{R}$ in three steps, first defining two
approximations $F,G$ and then gluing them together to get $\tilde{R}$.

Fix some $\gamma>0$, very small. Let $\phi+\psi$ be a $C^\infty$
partition of unity on $S^1$ associated to the sets $\{\theta \in
[0,\gamma)\}$ and $\{\theta \in (\gamma-\eta, \pi/2]\}$ where
$\theta$ is the angle on the circle (for some very small
$0<\eta<\gamma$). We define $F(z) =\phi(z)R(z)
+\psi(z)R(e^{i\gamma})$ on $\{\theta\in [0,\pi/2]\}$ : $F$ is
equal to $R$ on $\{\theta\in [0,\gamma-\eta]\}$ and to
$R(e^{i\gamma})$ on $\{\theta\in [\gamma,\pi/2]\}$. In particular,
the spectrum of $F(z)$ is ``almost the same'' as the spectrum
of $R(1)$ if $\gamma$ is small enough.

We define in the same way $F$ on $\{\theta\in [-\pi/2,0]\}$, equal
to $R(e^{-i \gamma})$ on $\{\theta\in [-\pi/2,-\gamma]\}$ and to
$R$ on $\{\theta\in [-\gamma+\eta,0]\}$.

Finally, we construct $F$ on the remaining half-circle by
symmetrizing, i.e. $F(e^{i(\pi/2+a)})=F(e^{i(\pi/2-a)})$, to
ensure that everything fits well.

Provided $\gamma$ is small enough, there is a well defined eigenvalue
close to $1$ for every $F(z)$,
depending continuously on $z$, which we denote by $\rho(z)$. The
problem would be solved if $\rho(z)\not=1$ for $z\not=1$, which
is not the case since $\rho(-1)=\rho(1)=1$. Consequently, we have
to perturb $\rho$ a little. Let $\nu$ be a $C^\infty$ function
on $\{\theta\in [\pi/2,3\pi/2]\}$, arbitrarily close to
$\rho$ and which does not take the value $1$.
On $\{\theta\in
[\pi/2+\eta,3\pi/2-\eta] \}$, we define
$G(z)=\frac{\nu(z)}{\rho(z)}F(z)$ : its
  eigenvalue close to $1$ is $\nu(z)\not=1$. Finally, we glue
  $F$ and $G$ together on $\{\theta\in [\pi/2,\pi/2+\eta]\}$ and
$\{\theta\in [3\pi/2-\eta,3\pi/2]\}$ with a partition of unity, as
above. As the spectrum of $F(e^{i\pi/2})=R(e^{i\gamma})$ does not
contain $1$, the gluing will not give an eigenvalue equal to $1$
if we choose $\eta$ small enough and $\nu$ close enough to $\rho$.

Note finally that, since $\frac{R(z)-R(1)}{z-1}\in A(\boBB)$, we also
have $\frac{\tilde{R}(z)-\tilde{R}(1)}
{z-1}\in A(\boBB)$: to obtain $\tilde{R}(z)$, we have
multiplied by $C^\infty$ functions, which are in $A(\C)$.

\emph{Step 3: Writing $\tilde{S}(z)=\frac{I-\tilde{R}(z)}{1-z}$, then
$\tilde{S}(z)^{-1} \in A(\boBB)$.}

By construction, $\frac{\tilde{R}(z)-\tilde{R}(1)}{z-1} \in
A(\boBB)$, and $\tilde{R}(z)\in A(\boBB)$. Equation
\eqref{pz-p1_sommable} (with tildes everywhere) implies that
$\frac{\tilde{P}(z)-\tilde{P}(1)}{z-1} \in A(\boBB)$ (use first
Wiener's Lemma Theorem \ref{wiener_classique} on $uI-\tilde{R}(z)$, then
multiply and integrate, as in the proof of Theorem
\ref{wiener_general}). In particular, $\tilde{P}(z) \in A(\boBB)$.

Let $\xi$ be a bounded linear functional on $\boBB$ satisfying
$\xi(P(1))=1$ and $\forall
z\in S^1,\ \xi(\tilde{P}(z))\not =0$ (take $\gamma$ small
enough in the construction of $\tilde{R}$ to ensure this). Then
$\xi(\tilde{P}(z))\in A(\C)$, thus $1/\xi(\tilde{P}(z))\in
A(\C)$ by Wiener's Lemma \ref{wiener_classique}.
Applying $\xi$ to both sides of
Equation \eqref{quotient_tend_vers_mu} (with tildes) and dividing
by $\xi(\tilde{P}(z))$, we obtain that
$\frac{1-\tilde{\lambda}(z)}{1-z} \in A(\C)$. As this function is
everywhere non-vanishing, we can use once more Wiener's Lemma and
get that $\frac{1-z}{1-\tilde{\lambda}(z)} \in A(\C)$.

Finally, we use Equation \eqref{S-1_converge_en_1} (with tildes
everywhere). Wiener's Lemma ensures that
$(I-\tilde{R}(z)\tilde{Q}(z))^{-1} \in A(\boBB)$ since
$I-\tilde{R}(z)\tilde{Q}(z)$ is in $A(\boBB)$ and pointwise
invertible. Thus, the right hand side is in $A(\boBB)$, which
ends the proof.

\emph{Step 4: $S(z)^{-1}\in A(\boBB)$.}

Let $\phi+\psi$ be a $C^\infty$ partition of unity on $S^1$, such
that $\tilde{R}=R$ on the support of $\phi$. Then
  \begin{equation*}
  \left(\frac{I-R(z)}{1-z}\right)^{-1}
  =\phi(z) \left(\frac{I-\tilde{R}(z)}{1-z}\right)^{-1}
  +\psi(z) \left(\frac{I-R(z)}{1-z}\right)^{-1}.
  \end{equation*}

The first term on the right-hand side is in $A(\boBB)$ according to Step 3.
For the second term, we can modify $\frac{I-R(z)}{1-z}$ outside
of the support of $\psi$ so that it is everywhere defined and
invertible, using the same techniques as in the construction of
$\tilde{R}$. Wiener's Lemma gives that its inverse is in
$A(\boBB)$, which concludes the proof.

\emph{Step 5: conclusion.}

$I-R(z)$ is a power series that is invertible everywhere on
$\D$. According to \cite[Lemma VII.6.13]{dunford_schwartz:1}, its
pointwise inverse is also analytic on $\D$, whence it is
a power series, that can be written as
$\sum T_n z^n$. Multiplying by $1-z$, we get $S(z)^{-1}=\sum (T_n
-T_{n-1}) z^n$. Thus, if $r<1$, we have $T_n-T_{n-1}=\frac{1}{2\pi r^n} \int
e^{-in\theta} S(re^{i\theta})^{-1} \dd\theta$. As $S(z)^{-1}$ is
continuous on the whole disk $\overline{\D}$, we obtain by letting
$r$ tend to $1$ that $T_n-T_{n-1}$ is the $n^{\text{th}}$ Fourier
coefficient of $S(z)^{-1}$ on $S^1$. But we already know that
$S(z)^{-1} \in A(\boBB)$, which proves that $\sum
\norm{T_n-T_{n-1}} <\infty$.
\end{proof}

\begin{rmq}
Theorem \ref{sarig_first_main}
is in fact the main technical ingredient in the proof
of Sarig in \cite{sarig:decay}, but he had to assume $a_n\sim
\frac{1}{n^\beta}$ with $\beta>2$.
Thus, it is possible using
Theorem \ref{sarig_first_main} to extend Sarig's results on lower
bounds for the decay of correlations (Theorem 2 and Corollaries 1 and
2 in \cite{sarig:decay}), even when the mass of the set of
points returning at time $n$ is not of the order $1/n^\beta$. For
example, it is not hard to generalize his results to maps of the
interval with a very neutral fixed point, such as $x(1+x \log^2
x)$. With a little work, this shows that all the upper bounds on
the decay of correlations obtained by Holland in \cite{holland}
are in fact optimal.
\end{rmq}

The following lemma will be needed later on in the proof of Theorem
\ref{TCLabstrait}.
\begin{lem}
\label{A_tend_0} Under the same hypotheses as in Theorem
\ref{sarig_first_main}, we have
$\left(\frac{I-R(z)}{1-z}\right)^{-1}=\frac{1}{\mu}P+(1-z)A(z)$,
where $A(z)=\sum A_n z^n$ with $A_n\to 0$ as $n\to \infty$.
\end{lem}
\begin{proof}
We use the method of Sarig in \cite{sarig:decay}, starting from
our analogue of his first main lemma, Theorem \ref{sarig_first_main}.
More precisely, Sarig constructs
a polynomial $R_B(z)$, approximating $R(z)$, such that
$\frac{R(1)-R_B}{1-z}$ and $\frac{1}{1-z}
\left(\frac{R(1)-R_B}{1-z}-R'(1)\right)$ are polynomials.
Moreover, $I-R_B(z)$ is invertible for $z\in \D$, and
$\left(\frac{I-R_B(z)}{1-z}\right)^{-1}=\frac{1}{\mu}P
+(1-z)B(z)$ where $B(z)=\sum B_n z^n$ and $B_n =o(\kappa^n)$ for
some $\kappa<1$. The proof of this result (\cite[Lemma
6]{sarig:decay}) relies only on his first main lemma, which we
have already proved in our context.

Then, he uses the following perturbative expansion of
$S(z)^{-1}=\left(\frac{I-R(z)}{1-z}\right)^{-1}$ around
$S_B(z)^{-1}=\left(\frac{I-R_B(z)}{1-z}\right)^{-1}$:
  \begin{equation*}
  S(z)^{-1}
  =S_B(z)^{-1}+S_B(z)^{-1} (S_B(z)-S(z)) S(z)^{-1}.
  \end{equation*}
To prove the lemma, we only have to show that $S_B(z)^{-1}
(S_B(z)-S(z)) S(z)^{-1}$ can be written as $(1-z) C(z)$ for some
$C(z)=\sum C_n z^n$ with $C_n \to 0$ (indeed, $A_n=B_n+C_n$).

We first study $\frac{S_B-S}{1-z}$. We have
  \begin{equation*}
  S_B-S=\left(\frac{R(1)-R_B}{1-z}-R'(1)\right)
  +\left(R'(1)-\frac{R(1)-R}{1-z}\right)
  =I+II.
  \end{equation*}
As $I/(1-z)$ is a polynomial (by construction of $R_B$), we can
forget this term. For the other term,
  \begin{equation*}
  II=\sum_{k=0}^\infty (1-z^k)\sum_{n=k+1}^\infty R_n.
  \end{equation*}
Thus, writing $D_n=\sum_{k=n+1}^\infty \sum_{l=k+1}^\infty R_l$,
we have $\frac{II}{1-z}=\sum D_n z^n$. The summability of
$a_k=\sum_{l=k+1}^\infty \norm{R_l}$ implies that $D_n\to 0$. We
have shown that
  \begin{equation*}
  \frac{S_B-S}{1-z}= \sum E_n z^n
  \end{equation*}
where $E_n \to 0$.

We already know that $S_B(z)^{-1}$ and $S(z)^{-1}$ have summable
coefficients. Thus, the lemma will be proved if we show that, if
$E_n \to 0$ and $F_n$ is summable, then the coefficient $(EF)_n$
of $z^n$ in  $\left(\sum E_n z^n\right)\left(\sum F_n z^n\right)$
tends to $0$ as $n\to\infty$.

We have $(EF)_n=\sum_{k=0}^n E_k F_{n-k}$. We fix $\epsilon>0$,
and $N$ such that if $k\geq N$, $\norm{E_k}\leq \epsilon$. If $K$
is greater than $\norm{E_k}$ and $\sum \norm{F_k}$, we get for
$n\geq N$
  \begin{equation*}
  \norm{(EF)_n}\leq \sum_{k=0}^{N-1} \norm{E_k}\norm{F_{n-k}} +\sum_{k=N}^n
  \epsilon \norm{F_{n-k}}
  \leq K \sum_{l=n-N}^\infty \norm{F_l} + K\epsilon.
  \end{equation*}
For $n$ large enough, this is less than $2K\epsilon$.
\end{proof}

\subsection{Invertibility of $I-R(z,t)$}

From this point on, and until the end of the proof of Theorem
\ref{TCLabstrait}, the notations and
hypotheses of all the results will be those of
Theorem \ref{TCLabstrait}.

In this section, we prove the invertibility of $I-R(z,t)$ when
$(z,t)\not =(1,0)$ and $t$ is small enough. The problem occurs for
$z$ close to $1$, where $R(z,t)$ has a well defined eigenvalue
$\lambda(z,t)$ close to $1$. Let us first study this eigenvalue.

\begin{lem}
\label{ecriture_valeur_propre} In a neighborhood of $(1,0)$, it
is possible to write
$\lambda(z,t)=1+(z-1)(\mu+a(z,t))-(c+b(t))M(|t|)$ where the
functions $a$ and $b$ tend to $0$ respectively at $(1,0)$ and $0$.
\end{lem}
\begin{proof}
As
  \begin{equation*}
  \frac{R(z,t)-R(1,t)}{z-1}
  =\sum_{n=0}^\infty z^n \left( \sum_{k>n}R_k(t)\right),
  \end{equation*}
the summability assumption implies that $(z,t)\mapsto
\frac{R(z,t)-R(1,t)}{z-1}$ is continuous at $t=0$. Writing
  \begin{equation*}
  \frac{P(z,t)-P(1,t)}{z-1}=\frac{1}{2i\pi}\int_{|u-1|=\delta} \frac{1}{uI-R(z,t)}
  \frac{R(z,t)-R(1,t)}{z-1}\frac{1}{uI-R(1,t)}\dd u,
  \end{equation*}
we also get the continuity of
 $(z,t)\mapsto \frac{P(z,t)-P(1,t)}{z-1}$ at $t=0$.
Thus, applying a bounded linear functional to the equation
  \begin{equation*}
  \frac{\lambda(z,t)-\lambda(1,t)}{z-1}P(z,t)
  =\frac{R(z,t)-R(1,t)}{z-1}P(z,t)+(R(1,t)-\lambda(1,t))
  \frac{P(z,t)-P(1,t)}{z-1},
  \end{equation*}
we get also the continuity of
$F(z,t)=\frac{\lambda(z,t)-\lambda(1,t)}{z-1}$. Moreover,
$F(z,0)\to \mu$ when $z\to 1$ (Equation \eqref{converge_vers_mu}),
which makes it possible to write
$F(z,t)=\mu+a(z,t)$ with $a(z,t)$ vanishing at $(1,0)$.
Moreover, by assumption, $\lambda(1,t)=1-(c+b(t))M(|t|)$ where $b$ vanishes
at $0$. Hence,
  \begin{equation*}
  \lambda(z,t)=\lambda(1,t)+(z-1)F(z,t)
  =1-(c+b(t))M(|t|) +(z-1)(\mu+a(z,t)).
  \end{equation*}
\end{proof}

\begin{prop}
\label{inversible_tout_D} There exists $\epsilon_0>0$ such that
$\forall (z,t)\in (\overline{\D}×
[-\epsilon_0,\epsilon_0])-\{(1,0)\}$, $I-R(z,t)$ is invertible.
\end{prop}
\begin{proof}
Assume that there is a neighborhood $V$ of $(1,0)$ such that, in
$V-\{(1,0)\}$,
$I-R(z,t)$ is invertible.
We can assume that $V=U× [-\epsilon_1,\epsilon_1]$. For
$w\in \overline{\D}-U$, $I-R(w)$ is invertible, thus $I-R(z,t)$
is invertible in a small neighborhood of $(w,0)$, which can be
taken of the form $U_w × [-\epsilon_w,\epsilon_w]$. By compactness,
$\overline{\D}-U$ is covered by a finite number of the $U_w$, say
$U_{w_1},\ldots,U_{w_n}$. To conclude, just take
$\epsilon_0=\min(\epsilon_1,
\epsilon_{w_1},\ldots,\epsilon_{w_n})$.

So, we just have to construct
a neighborhood $V$ as above.
We have to show that
the eigenvalue $\lambda(z,t)$ of $R(z,t)$ is different from $1$
when $(z,t)$ is close to but different from $(1,0)$. This
essentially comes from the asymptotic expansion of $\lambda(z,t)$
given by Lemma \ref{ecriture_valeur_propre}, but we will write it
more carefully.

Let $K$ be such that $\frac{\mu}{4}>\frac{3|c|}{2K}$. The
neighborhood $V$ will be the set of $(z,t)$ such
that $|b(t)|\leq \Re(c)/2$ and
$|a(z,t)|\leq\min\left(\frac{\mu}{4},\frac{\Re(c)}{2K}\right)$.

We fix $(z,t)\in V\cap (\overline{\D}× \R)$. Assuming that
$\lambda(z,t)=1$, we show that $(z,t)=(1,0)$. There are two
possibilities, namely $|\Im(z-1)|\geq |z-1|/2$ or $|\Re(z-1)|\geq
|z-1|/2$. The latter case, being easier, is left to the reader. We can
assume for example that $\Im(z-1)\geq |z-1|/2$.

As $\Re(\mu(z-1))\leq 0$ since $\mu>0$ and $z\in \overline{\D}$,
we have
  \begin{align*}
  1=\Re(\lambda&(z,t))
  =1+\Re(\mu(z-1))+\Re(a(z,t)(z-1))
  -(\Re(c)+\Re(b(t))M(|t|)
  \\&
  \leq
  1+|a(z,t)||z-1|-(\Re(c)-|b(t)|) M(|t|)
  \leq
  1+\frac{\Re(c)}{2K} |z-1| -\frac{\Re(c)}{2} M(|t|).
  \end{align*}
Thus, $|z-1|\geq K M(|t|)$. The imaginary part then satisfies
  \begin{align*}
  0=\Im(\lambda(z,t))&
  \geq (\mu-|a(z,t)|) \Im(z-1)-(|c|+|b(t)|)M(|t|)
  \geq \frac{\mu}{2} \frac{|z-1|}{2}- \frac{3|c|}{2} \frac{|z-1|}{K}.
  \end{align*}
The definition of $K$ implies that $|z-1|=0$, which in turn gives
$t=0$ since $|z-1|\geq K M(|t|)$ and $M$ vanishes only at $0$.
%
\end{proof}

In particular, for $|t|<\epsilon_0$, $I-R(z,t)$ is a power series
 invertible everywhere on $\D$, whence its inverse can be written as
$\sum T_{n,t}z^n$ (see Step 5 of the proof of Theorem
\ref{sarig_first_main}).
Moreover, for $t\not= 0$, everything can be
continuously extended to $\partial \D$, which implies that the
$T_{n,t}$ are the Fourier coefficients of the function
$(I-R(z,t))^{-1}$ (considered as a function from $S^1$ to
$\boBB$).

\subsection{Convergence of $\tilde{R}(z,t)$}

In this section, we will work exclusively on $S^1$
($z$ will thus denote a point in $S^1 \subset \C$
and \emph{not} in $\D$).

To apply Fourier series methods to $R(z,t)$, we will need it to
have a well defined eigenvalue close to $1$, which is not \emph{a
priori} true outside of a neighborhood of $(1,0)$.

{\bf Construction of $\tilde{R}(z,t)$.}

We use the function $\tilde{R}(z)$ that has been constructed in
Step 2 of the proof of Theorem \ref{sarig_first_main}. We set
$\tilde{R}(z,t)=\tilde{R}(z) -\left( R(z)-\sum_{n=1}^\infty z^n
R_n(t)\right)$. The map $t \mapsto \tilde{R}(\cdot,t)$ is then continuous
at $t=0$
(for the $A(\boBB)$-norm, whence for the supremum norm), and
$\tilde{R}(z,t)$ has a well defined eigenvalue
$\tilde{\lambda}(z,t)$ close to $1$. Moreover, for $z$ in a
neighborhood of $1$, $\tilde{\lambda}(z,t)=\lambda(z,t)$, which
is different from $1$ when $(z,t)\not =(1,0)$, and outside of
this neighborhood $\tilde{\lambda}(z,t)$ is close to
$\tilde{\lambda}(z,0)$ and is thus different from $1$ if $t$ is
small enough.

The main result of this section will be the following proposition
(recall that the hypotheses are those of Theorem \ref{TCLabstrait}):

\begin{prop}
\label{prop_conv_tildeR}
We have
  \begin{equation*}  
  \norm{\left(\frac{I-\tilde{R}(z,t)}{1-(1-\frac{c}{\mu}M(|t|))z}\right)^{-1} -
  \left( \frac{I-\tilde{R}(z)}{1-z}\right)^{-1}}_{A(\boBB)}
  \xrightarrow[t\to 0]{} 0.
  \end{equation*}
\end{prop}
\begin{proof}
Write $\tilde{P}(z,t)$ for the spectral projection
associated to the eigenvalue $\tilde{\lambda}(z,t)$, and
$\tilde{Q}(z,t)=I-\tilde{P}(z,t)$. We then have
  \begin{equation*}
  I-\tilde{R}(z,t)=(1-\tilde{\lambda}(z,t))\tilde{P}(z,t)+
  (I-\tilde{R}(z,t)\tilde{Q}(z,t))\tilde{Q}(z,t).
  \end{equation*}
Thus,
  \begin{multline*}
  \left(\frac{I-\tilde{R}(z,t)}{1-(1-\frac{c}{\mu}M(|t|))z}\right)^{-1}
  \\
  =\frac{1-(1-\frac{c}{\mu}M(|t|))z}{1-\tilde{\lambda}(z,t)}\tilde{P}(z,t)
  +\left(1-(1-\frac{c}{\mu}M(|t|))z\right)(I-\tilde{R}(z,t)\tilde{Q}(z,t))^{-1}
  \tilde{Q}(z,t).
  \end{multline*}
As $\tilde{R}(z,t)\to \tilde{R}(z)$ as $t\to 0$, Wiener's Lemma
Theorem
\ref{wiener_general} gives the convergence of $\tilde{P}(z,t)$
and $\tilde{Q}(z,t)$ respectively to $\tilde{P}(z)$ and
$\tilde{Q}(z)$ (in the sense of Fourier coefficients). 
We also obtain (again by Wiener's Lemma) the
convergence of $(I-\tilde{R}(z,t)\tilde{Q}(z,t))^{-1}$ to
$(I-\tilde{R}(z)\tilde{Q}(z))^{-1}$, since
$I-\tilde{R}(z)\tilde{Q}(z)$ is
invertible for every $z\in S^1$.


It only remains to prove
that
$\frac{1-(1-\frac{c}{\mu}M(|t|))z}{1-\tilde{\lambda}(z,t)} \to
\frac{1-z}{1 -\tilde{\lambda}(z)}$ (in $A(\C)$),
which
is given by the following lemma and Wiener's Lemma (since
inversion is continuous in the Banach algebra $A(\C)$).
\end{proof}

\begin{lem}
\label{lemme_crucial_estime_vp}
We have
  \begin{equation*}
  \norm{\frac{1-\tilde{\lambda}(z,t)}
  {1-(1-\frac{c}{\mu}M(|t|))z}-\frac{1-\tilde{\lambda}(z)}{1-z}
  }_{A(\C)}\xrightarrow[t\to 0]{} 0.
  \end{equation*}
\end{lem}
\begin{proof}
We will usually work with the circle $S^1$ embedded in $\C$,
but it will be
convenient to consider the circle as $\T=\R/2\pi\Z$, in
which case the variable will be denoted by $\theta$.

We recall some results on Fourier series. For
$|\theta|\leq \pi$ and small $\epsilon>0$, write
$\Delta_\epsilon(\theta)=\sup\left(0,1-\frac{|\theta|}
{\epsilon}\right)$: this function is supported in
$\{\theta\in [-\epsilon,\epsilon]\}$. We write also
$V_\epsilon=2\Delta_{2\epsilon}-\Delta_\epsilon$: this function
is equal to $1$ on $\{\theta\in[-\epsilon,\epsilon]\}$ and vanishes outside of
$\{\theta\in[-2\epsilon,2\epsilon]\}$. We will use $V_\epsilon$ and
$1-V_\epsilon$ as a $C^0$ partition of unity. By \cite[page
56]{kahane:sommable}, we have $\norm{V_\epsilon}_{A(\C)} \leq 3$. Note
that $V_\epsilon$ is piecewise $C^1$, its derivative being bounded by
$1/\epsilon$.

If $g:\T\to \C$ is continuous and piecewise $C^1$, then (\cite[page 56,
Equation (1)]{kahane:sommable})
  \begin{equation}
  \label{majore_somme_par_derivee}
  \norm{g}_{A(\C)} \leq |c_0(g)|+\left(\frac{\pi}{12}\int_{-\pi}^{\pi}
    |g'(\theta)|^2\dd t\right)^{1/2}.
  \end{equation}

Finally, if $f:S^1\to \C$ is in  $A(\C)$ and satisfies $f(1)=0$,
then (\cite[page 56]{kahane:sommable})
  \begin{equation}
  \label{lemme_kahane}
  \norm{fV_\epsilon}_{A(\C)}\xrightarrow[\epsilon \to 0]{}0.
  \end{equation}

Fix $|t|<\epsilon_0$ and write $\psi_t=V_{G(|t|)}$ for
 $G(t)=M(t)^{1/3}$, and
$\phi_t=1-\psi_t$, on $S^1$ or $\T$. Setting
$A_t(z)=\frac{1-\tilde{\lambda}(z,t)}{1-(1-\frac{c}{\mu}M(|t|))z}$, we
have to check that $A_t(z)$ tends (in $A(\C)$) to
$B(z)=\frac{1-\tilde{\lambda}(z)}{1-z}$. We will prove that
$A_t(z)\psi_t(z) - B(z) \psi_t(z)\to 0$, and $A_t(z)\phi_t(z) -
B(z) \phi_t(z)\to 0$.

In this proof, we denote by $C(z)$ a function whose norm in
$A(\C)$ is bounded by a constant. For example, as
$\norm{V_t}_{A(\C)}\leq 3$, we sometimes replace $V_t(z)$ by
$C(z)$. As $\frac{1-z}{1-(1-\frac{c}{\mu}M(|t|))z} =
1-\frac{c}{\mu}M(|t|) \sum_{n=1}^\infty
(1-\frac{c}{\mu}M(|t|))^{n-1}z^n $, we also have
$\norm{\frac{1-z}{1-(1-\frac{c}{\mu}M(|t|))z}}_{A(\C)} \leq
1+\frac{
\frac{|c|}{\mu}M(|t|)}{1-\left|1-\frac{c}{\mu}M(|t|)\right|}$,
which is bounded when $t\to 0$ (since $\Re c>0$). Thus, we
may also replace $\frac{1-z}{1-(1-\frac{c}{\mu}M(|t|))z}$ by
$C(z)$.

Using exactly the same equations as in the proof of Lemma
\ref{ecriture_valeur_propre} (with tildes everywhere) and
Wiener's Lemma, we prove that
$K(z,t)=\frac{\tilde{\lambda}(z,t)-\tilde{\lambda}(1,t)}{z-1}
-\frac{\tilde{\lambda}(z)-\tilde{\lambda}(1)}{z-1}$ tends to
$0$ in $A(\C)$ when $t\to 0$.

\emph{Step 1 (close to $z=1$): $A_t(z)\psi_t(z) -B(z) \psi_t(z)$ tends to $0$
in $A(\C)$-norm.}

We proved in the third step of the proof of Theorem
\ref{sarig_first_main} that
$\frac{\tilde{\lambda}(z)-\tilde{\lambda}(1)}{z-1}\in A(\C)$. As this
quotient
tends to $\mu$ when $z\to 1$, we can write
$\frac{\tilde{\lambda}(z)-\tilde{\lambda}(1)}{z-1}=\mu+F(z)$ where
$F \in A(\C)$ and $F(1)=0$. We also write
$\tilde{\lambda}(1,t)=1-(c+b(t))M(|t|)$ (since
$\tilde{\lambda}(1,t)=\lambda(1,t)$, as $\tilde{R}(z)=R(z)$ in a
neighborhood of $1$). This gives
  \begin{equation*}
  \tilde{\lambda}(z,t)
    =1-(c+b(t))M(|t|)+(z-1)\left(K(z,t)+F(z)+\mu\right).
  \end{equation*}
Thus,
  \begin{multline*}
  A_t(z)=
  \frac{1-\tilde{\lambda}(z,t)}{1-(1-\frac{c}{\mu}M(|t|))z}
  \\
  =\mu+\frac{1-z}{1-(1-\frac{c}{\mu}M(|t|))z}\left(F(z)+K(z,t)+cM(|t|)\right)
  +\frac{M(|t|)}{1-(1-\frac{c}{\mu}M(|t|))z} b(t).
  \end{multline*}
Therefore,
  \begin{multline*}
  A_t(z)\psi_t(z)
  \\=\mu \psi_t(z)+C(z)F(z)\psi_t(z)+C(z)(K(z,t)+M(|t|))+
  \frac{M(|t|)}{1-(1-\frac{c}{\mu}M(|t|))z} b(t)\psi_t(z).
  \end{multline*}
A small computation shows that
$\frac{M(|t|)}{1-(1-\frac{c}{\mu}M(|t|))z}$ has a bounded norm in
$A(\C)$. Moreover, $\mu=B(z)-F(z)$, and the term $F(z) \psi_t(z)$
can be included in $C(z)F(z)\psi_t(z)$. Finally,
  \begin{equation*}
  A_t(z)\psi_t(z)-B(z)\psi_t(z)=C(z)F(z)\psi_t(z)+C(z)(M(|t|)+K(z,t)+b(t)).
  \end{equation*}
In the rightmost term, everything tends to $0$ as $t\to 0$.
Moreover, Equation \eqref{lemme_kahane} ensures that
$F(z)\psi_t(z) \to 0$. This concludes the first step.

\emph{Step 2 (away from $z=1$): $A_t(z)\phi_t(z) -B(z) \phi_t(z)$ tends to $0$
in $A(\C)$-norm.}

We start from
  \begin{multline*}
  \frac{1-\tilde{\lambda}(z,t)}{1-(1-\frac{c}{\mu}M(|t|))z}
  -\frac{1-\tilde{\lambda}(z)}{1-z}
  \\
  =\frac{1-\tilde{\lambda}(z)}{1-z}\left(\frac{1-z}
     {1-(1-\frac{c}{\mu}M(|t|))z}-1 \right)
  +\frac{1-z}{1-(1-\frac{c}{\mu}M(|t|))z}
   \frac{\tilde{\lambda}(z)-\tilde{\lambda}(z,t)}{1-z}.
  \end{multline*}
We will see that each term on the right, multiplied by $\phi_t$,
tends to $0$ in $A(\C)$ when $t\to 0$.

We will use that, on the support of $\phi_t$, we have $|z-1| \geq
C G(|t|)$, which implies that
$|1-(1-\frac{c}{\mu}M(|t|))z|\geq C G(|t|)-D M(|t|) \geq C
  G(|t|)$ (for a smaller constant $C$) since $M(t)=o(G(t))$. In the
following, $C$ will denote a generic constant.

\emph{The first term.}

As $\frac{1-\tilde{\lambda}(z)}{1-z} \in A(\C)$, it is enough to
prove the convergence of $f=\left(\frac{1-z}
     {1-(1-\frac{c}{\mu}M(|t|))z}-1 \right)\phi_t(z)$ to $0$.
We will use Equation \eqref{majore_somme_par_derivee}.

We have
  \begin{equation*}
  |f(z)|\leq \frac{\frac{|c|}{\mu}M(|t|)}{C G(t)}\leq C
  M(|t|)^{2/3}.
  \end{equation*}
In particular, $|c_0(f)|\leq C M(|t|)^{1/3}$.

For the derivative, writing $f(\theta)=h(\theta)\phi_t(\theta)$,
we have $f'=h'\phi_t+h \phi_t'$. As $|h| \leq CM(|t|)^{2/3}$ and
$|\phi_t'|\leq 1/G(|t|)$, we get $|h\phi_t'|\leq CM(|t|)^{1/3}$.

To control $h'$, we write
  \begin{equation*}
  h(\theta)=-\frac{\frac{c}{\mu}M(|t|)}{1-\frac{c}{\mu}M(|t|)}+
  \frac{\frac{c}{\mu}M(|t|)}{1-\frac{c}{\mu}M(|t|)}
  \frac{1}{1-(1-\frac{c}{\mu}M(|t|))e^{i\theta}}.
  \end{equation*}
Thus,
  \begin{equation*}
  h'(\theta)=i\frac{c M(|t|)}{\mu}
  \frac{1}{(1-(1-\frac{c}{\mu}M(|t|))e^{i\theta})^2}.
  \end{equation*}
On the support of $\phi_t$, the modulus of the
denominator on the right is $\geq
(CG(|t|))^2$, which gives $|h'|\leq C M(|t|)^{1/3}$.

We have shown that $|f'|\leq C M(|t|)^{1/3}$. As $|c_0(f)|\leq
C M(|t|)^{1/3}$, Equation \eqref{majore_somme_par_derivee} gives that
$\norm{f}_{A(\C)}\leq C M(|t|)^{1/3}$.

\emph{The second term.}

As $\frac{1-z}{1-(1-\frac{c}{\mu}M(|t|))z}$ is bounded in $A(\C)$,
it is enough to show that
$g(z)=\frac{\tilde{\lambda}(z,t)-\tilde{\lambda}(z)}{1-z}\phi_t(z)$
tends to $0$ in $A(\C)$.

Writing again
$K(z,t)=\frac{\tilde{\lambda}(z,t)-\tilde{\lambda}(1,t)}{z-1}
-\frac{\tilde{\lambda}(z)-\tilde{\lambda}(1)}{z-1}$ (this difference
tends to $0$ when $t\to 0$), we have
  \begin{equation*}
  g(z)=K(z,t)\phi_t(z)+(\tilde{\lambda}(1)
-\tilde{\lambda}(1,t))\frac{\phi_t(z)}{1-z}.
  \end{equation*}
The first term $K(z,t)\phi_t(z)$ tends to $0$. For the second
term, if $h(z)=\frac{\phi_t(z)}{1-z}$, we have $|c_0(h)|\leq
\frac{1}{C G(|t|)}$ and
  \begin{equation*}
  |h'(\theta)|=\left|i\frac{\phi_t(\theta)}{(1-e^{i\theta})^2}
+\frac{1}{1-e^{i\theta}} \phi_t'(\theta)\right|
  \leq
  \frac{1}{C^2G(|t|)^2}+\frac{1}{CG(|t|)}\frac{1}{G(|t|)}.
  \end{equation*}
Equation \eqref{majore_somme_par_derivee} gives that
$\norm{h}_{A(\C)} \leq \frac{C}{G(|t|)^2}=O(M(|t|)^{-2/3})$. As
$\tilde{\lambda}(1)-\tilde{\lambda}(1,t)=O(M(|t|))$, we finally
get $\norm{(\tilde{\lambda}(1)-\tilde{\lambda}(1,t))
  \frac{\phi_t(z)}{1-z}}_{A(\C)}=O(M(|t|)^{1/3})$, which concludes the
proof of Step 2 of the proof of Lemma \ref{lemme_crucial_estime_vp}.
\end{proof}

\subsection{Conclusion of the proof of Theorem \ref{TCLabstrait}}

\begin{prop}
\label{prop_convergence_R} We have
  \begin{equation*}
  \norm{\left(\frac{I-R(z,t)}{1-(1-\frac{c}{\mu}M(|t|))z}\right)^{-1} -
    \left( \frac{I-R(z)}{1-z}\right)^{-1}}_{A(\boBB)}
    \xrightarrow[t\to 0]{} 0.
  \end{equation*}
\end{prop}
\begin{proof}
Let $\chi_1,\chi_2$ be a partition of unity such that
$R=\tilde{R}$ on the support of $\chi_1$, and set
$A(z,t)=\left(\frac{I-R(z,t)}{1-(1-\frac{c}{\mu}M(|t|))z}\right)^{-1}$.

Proposition \ref{prop_conv_tildeR} shows that $\chi_1(z) A(z,t)
\to \chi_1(z) A(z,0)$. For the term $\chi_2(z) A(z,t)$, we can
safely modify $R$ on a small neighborhood of $1$ so that $I-R(z)$
is everywhere invertible (this does not change $\chi_2(z)
A(z,t)$). Wiener's Lemma and the continuity of inversion then give
that $\chi_2(z) A(z,t)\to
\chi_2(z) A(z,0)$.
\end{proof}

\begin{proof}[Proof of Theorem \ref{TCLabstrait}]
We now conclude the proof of Theorem \ref{TCLabstrait}, i.e. we
construct functions $\epsilon(t)$ and $\delta(n)$ satisfying Equation
\eqref{EquationTCLabstrait}.

Proposition \ref{prop_convergence_R} shows that it is possible to
write
  \begin{equation*}
  \left(\frac{I-R(z,t)}{1-(1-\frac{c}{\mu}M(|t|))z}\right)^{-1}=
    \left( \frac{I-R(z)}{1-z}\right)^{-1}+F_t(z)
  \end{equation*}
where $F_t(z)\in A(\boBB)$ tends to $0$ when $t\to 0$. Moreover,
Lemma \ref{A_tend_0} gives that
  \begin{equation*}
  \left(\frac{I-R(z)}{1-z}\right)^{-1}=\frac{1}{\mu}P+(1-z)A(z)
  \end{equation*}
where $A(z)=\sum A_n z^n$ with $A_n\to 0$. We get
  \begin{multline*}
  \sum T_{n,t}z^n=(I-R(z,t))^{-1}
  \\
  =\frac{1}{1-(1-\frac{c}{\mu}M(|t|))z} \frac{1}{\mu}P
  +\frac{1-z}{1-(1-\frac{c}{\mu}M(|t|))z}A(z)
  +\frac{1}{1-(1-\frac{c}{\mu}M(|t|))z}F_t(z).
  \end{multline*}
In the first term, the coefficient of $z^n$ is
$\frac{1}{\mu}\left(1-\frac{c}{\mu}M(|t|)\right)^n P$ and
corresponds to \eqref{EquationTCLabstrait}. We have to check that
the other terms are smaller than $\epsilon(t)+\delta(n)$.

Let us denote by $\delta(n)$ the maximum for $t\in
[-\epsilon_0,\epsilon_0]$ of the coefficient $c_n(t)$ of $z^n$ in
$\frac{1-z}{1-(1-\frac{c}{\mu}M(|t|))z}A(z)$. We will prove that
$\delta(n)\to 0$ as $n\to\infty$.
We have
  \begin{equation*}
  \frac{1-z}{1-(1-\frac{c}{\mu}M(|t|))z}=1-\frac{c}{\mu}
   M(|t|)\sum_{n=1}^\infty
  \left(1-\frac{c}{\mu}M(|t|)\right)^{n-1}z^n.
  \end{equation*}
Thus,
  \begin{equation*}
  c_n(t)=A_n-\frac{c}{\mu}M(|t|)\sum_{k=0}^{n-1}\left(1-\frac{c}{\mu}
   M(|t|)\right)^{n-k-1}A_k.
  \end{equation*}
Note that, since $\Re(c)>0$ and $M(t)\to 0$ when $t\to 0$, we have
 $\left|1-\frac{c}{\mu}M(|t|)\right|\leq 1-\frac{\Re
  c}{2\mu}M(|t|)$ when $t$ is small enough.
Let $C$ be greater than the $\norm{A_n}, n\in \N$. We fix
$\epsilon>0$, and $N$ such that $\forall k\geq N$,
$\norm{A_k}\leq \epsilon$. Then, if $n\geq N$,
  \begin{align*}
  \norm{c_n(t)}&\leq \norm{A_n}+\frac{|c|}{\mu}
  M(|t|)\sum_{k=0}^{N-1}\left(1-\frac{\Re c}{2\mu}M(|t|)\right)^{n-k-1}
  \norm{A_k}
  \\&\hphantom{\leq \;}+\frac{|c|}{\mu}M(|t|)\sum_{k=N}^{n-1} \left(1-
  \frac{\Re c}{2\mu} M(|t|)\right)^{n-k-1}\epsilon
  \\&
  \leq \epsilon +N\frac{|c|}{\mu}
  CM(|t|) \left(1-\frac{\Re c}{2\mu}M(|t|)\right)^{n-N} +E\epsilon.
  \end{align*}
As $M(|t|)\left(1-\frac{\Re c}{2\mu}M(|t|)\right)^{n-N}$
tends uniformly to $0$ on
$[-\epsilon_0,\epsilon_0]$ (if $\epsilon_0$
is small enough) when $n\to\infty$, we get that for $n$ large
enough, $\norm{c_n(t)}\leq (2+E)\epsilon$, and $\delta(n)\leq
(2+E)\epsilon$.

Let us finally bound the coefficients of
$\frac{1}{1-(1-\frac{c}{\mu}M(|t|))z}F_t(z)$, which will give the
$\epsilon(t)$.  The Fourier coefficients of
$\frac{1}{1-(1-\frac{c}{\mu}M(|t|))z}$ are $\leq 1$, which implies
that the coefficient of $z^n$ in the product is at most
$\sum_{k=0}^n \norm{F_t(z)_k} \leq \norm{F_t}_{A(\boBB)}$. Since
$\norm{F_t}_{A(\boBB)} \to 0$ when $t\to 0$, we can take
$\epsilon(t)=\norm{F_t}_{A(\boBB)}$ to conclude.
\end{proof}

\section{Results in the uniform case}
\label{section:uniform_case}

In Sections \ref{section_quelques_lemmes} and
\ref{section_lemmes_perturbatifs}, we state some technical estimates
on Markov maps, and we use them in Sections
\ref{section_TCL_classique} and \ref{section_stable_classique} to
derive an expansion of the eigenvalue of a perturbed transfer
operator, which will subsequently be used in the proofs of Theorems
\ref{TCL_markov} and
\ref{thm_loi_stable_markov}.

We work in the setting of Markov maps of Section
\ref{section_def_markov_maps}, and we will use freely the definitions
introduced here.

\subsection{Some technical estimates}
\label{section_quelques_lemmes}

In this section, $(Z,\boB,T,m,\delta)$ will be a Markov map
preserving a measure $m$ of finite mass. Write $t(x,y)$ for the
separation time of the points $x$ and $y$ (under the iteration of
$T$) and assume that the distortion of $T$ is $\tau$-Hölder for
some $\tau<1$. For $f:Z \to \C$, we will denote by $D_\tau f$ the
least $\tau$-Hölder constant of $f$ for points in the same
elements of the partition, i.e. $D_\tau f=\inf\{ c \tq \forall
x,y\in Z, t(x,y)\geq 1 \Rightarrow |f(x)-f(y)|\leq c
\tau^{t(x,y)}\}$. Let also $\boL_\tau$ be the space of bounded
functions on $Z$ with $D_\tau f<\infty$. It is a Banach space for
the norm $\norm{f}_{\boL_\tau}=\norm{f}_\infty +D_\tau f$.

We will need a property slightly weaker than the big image
property from Section \ref{section_def_markov_maps}.
We say that $Z$ is a \emph{bounded tower} if there exists
$q\in \N$, a subset $\delta_0\subset \delta$ (corresponding to the
basis of the tower) and a function $R: \delta_0 \to \{0,\ldots,q-1\}$,
such that $Z$ is isomorphic to the set
$\{(x,i) \tq x\in a\subset \delta_0, i<R(a)\}$, with the partition into
the sets $a \times \{i\}$ for $a\in \delta_0$ and $i<R(a)$. We require
$T$ to map $a \times \{i\}$ isomorphically onto $a\times\{i+1\}$ if
$i+1<R(a)$, while if $i+1=R(a)$ then $T(a\times\{i\})\subset
\bigcup_{d\in \delta_0}d$, i.e. the points get out of the tower
through the top and fall down to the basis. We will write
$\Delta_k=\{(x,k)\tq x\in a, R(a)>k\}$,
i.e. this is the set of points at height $k$ in the tower.


We will say that $T$ has the \emph{bounded tower big image property} if
$Z$ is a bounded tower, and the returns to the basis have measure
bounded away
from $0$, i.e. $\forall a\in \delta$, $Ta\subset \Delta_0
\Rightarrow m[Ta]>\eta$ for a constant $\eta>0$. Note that, when
$q=1$, i.e. there are no floors in the tower, we recover the usual
big image property.

In the rest of Section \ref{section_quelques_lemmes},
$T$ will have the bounded tower big image
property and $\log g_m$ will be $\tau$-Hölder.

\begin{lem}[Distortion Lemma]
\label{controle_distortion}
Writing $g_m^{(n)}(x)=g_m(x)\cdots g_m(T^{n-1}x)$ (this is the inverse
of the
jacobian of $T^n$), there exists a
constant $C$ such that $t(x,y)\geq n \Rightarrow \left|
  1-\frac{g_m^{(n)}(x)}{g_m^{(n)}(y)}\right| \leq C \tau^{t(x,y)-n}$.

In particular, there exists a constant $D$ such that,
if $\underline{d}=[d_0,\ldots,d_{n-1}]$ is a cylinder
of length $n$ and $x\in \underline{d}$, then $D^{-1} g_m^{(n)}(x) \leq
\frac{m[\underline{d}]}{m[T^n \underline{d}]} \leq D g_m^{(n)}(x)$.
\end{lem}
\begin{proof}
Let $x,y$ be such that $t(x,y)\geq n$. Then
  \begin{align*}
  \left|\log (g_m^{(n)}(x))-\log(g_m^{(n)}(y))\right|
  &\leq \sum_{i=0}^{n-1}\left| \log(g_m(T^ i x))-\log(g_m(T^i
  y))\right|
  \\&
  \leq \sum_{i=0}^{n-1} D_\tau(\log(g_m)) \tau^{t(x,y)-i}
  \leq \frac{D_\tau(\log(g_m))}{1-\tau} \tau^{t(x,y)-n+1}.
  \end{align*}
Taking exponentials, we get the first statement of the lemma.

In particular, if $\underline{d}$ is a cylinder of length $n$ and
$x,y\in \underline{d}$,
we have $\frac{g_m^{(n)}(x)}{g_m^{(n)}(y)} \leq D$ for some $D$.
Using that $1/g_m^{(n)}(y)$ it the jacobian of $T^n$, we integrate on
$y\in\underline{d}$ and get that $g_m^{(n)}(x)m[T^n \underline{d}]\leq D
m[\underline{d}]$.
The other inequality is analogous.
\end{proof}

If $\psi$ is a function on $Z$ and $\hat{T}$ denotes the transfer
operator associated to $T$, then $\hat{T}^n \psi(x)=\sum_{T^n(y)=x}
g_m^{(n)}(y)\psi(y)$. There is at most one preimage of $x$ in every
non-empty cylinder $\underline{d}=[d_0,\ldots,d_{n-1}]$, which we
denote by $d_0\ldots d_{n-1}x$. Set $M_{\underline{d}}\psi(x)=
g_m^{(n)}(d_0\ldots d_{n-1}x) \psi(d_0\ldots d_{n-1}x)$ if this point
is defined, and $0$ otherwise. Then $\hat{T}^n=\sum M_{\underline{d}}$
where the sum is over all non-empty cylinders $\underline{d}$ of
length $n$.

\begin{lem}
\label{lemme_sur_Mdbarre}
There exists a constant $B$ such that
$\forall \underline{d}=[d_0,\ldots,d_{n-1}]$, $\forall \psi: Z\to \C$,
  \begin{equation*}
  \norm{M_{\underline{d}}\psi}_{\boL_\tau} \leq B
  \frac{m[\underline{d}]}{m[T^n\underline{d}]}\left(
  \tau^n D_\tau\psi+\frac{1}
  {m[\underline{d}]} \int_{[\underline{d}]} |\psi|
      \dd m\right).
  \end{equation*}
\end{lem}
\begin{proof}
Let $x$ and $y$ be in the same partition element $a$. If $a$ is
not in the image of $T(d_{n-1})$, $(M_{\underline{d}}\psi)_{|a}=0$
and there is nothing to prove. Otherwise, write $x'=d_0\ldots
d_{n-1}x$, and $y'=d_0\ldots d_{n-1}y$.

First note that, for $z'\in \underline{d}$, we have $|\psi(y')-\psi(z')|\leq
D_\tau\psi\; \tau^n$,  which gives, after integrating on $z'$,
 that $|\psi(y')|\leq D_\tau\psi\; \tau^n
+\frac{1}{m[\underline{d}]}\int_{\underline{d}} |\psi|\dd m$. The
inequality $g_{m}^{(n)}(y') \leq D \frac{m[\underline{d}]}{m[T^n
\underline{d}]}$ from Lemma \ref{controle_distortion}
gives then the claimed bound on
$\norm{M_{\underline{d}} \psi}_\infty$.

Using
Lemma \ref{controle_distortion}, we have
  \begin{align*}
  |M_{\underline{d}}\psi(x)-M_{\underline{d}}\psi(y)|&
  \leq |g_{m}^{(n)}(x')| |\psi(x')-\psi(y')|+ |\psi(y')||g_{m}^{(n)}(y')|
  \left|1-\frac{g_{m}^{(n)}(x')}{g_{m}^{(n)}(y')}\right|
  \\&
  \leq D \frac{m[\underline{d}]}{m[T^n \underline{d}]}D_\tau\psi\;
   \tau^{t(x,y)+n}
   + |\psi(y')| D \frac{m[\underline{d}]}{m[T^n \underline{d}]}
  C \tau^{t(x,y)}
  \end{align*}
which gives the conclusion, using the aforementioned bound on $|\psi(y')|$.
\end{proof}

\begin{cor}
\label{big_donc_vpisolee} If $T$ has the bounded tower big image property
and is ergodic, then the associated transfer operator
$\hat{T}$ acts continuously on $\boL_\tau$,
and it has a simple isolated eigenvalue at $1$, the eigenspace
being the constant functions.
\end{cor}
\begin{proof}
The proof will use Lemma \ref{lemme_sur_Mdbarre}
and the fact that $\hat{T}^n=\sum
M_{\underline{d}}$, where the sum is over all cylinders of length $n$.

We first show the continuity of $\hat{T}$. Decompose
$\hat{T}=1_{\Delta_0} \hat{T} + (1-1_{\Delta_0}) \hat{T}$, where
the first term sees the points on the basis and the second one in
the floors. If $x$ is in the basis, the elements $a$ of the
partition such that $x\in Ta$ return to the basis, which implies
that they have a big image. We can thus sum the estimates given by
Lemma \ref{lemme_sur_Mdbarre} for cylinders of length $1$, and
forget the $m[T\underline{d}]$, to get that $\norm{1_{\Delta_0}
\hat{T} \psi}_{\boL_\tau}\leq B(\tau
\norm{\psi}_{\boL_\tau}+\norm{\psi}_1)$. On the other hand, if $x$
is not in the basis and $x'$ denotes the point just below, then $
\hat{T} \psi(x)=\psi(x')$. This implies that
$\norm{(1-1_{\Delta_0}) \hat{T}}_{\boL_\tau}\leq1$, and concludes
the proof of the continuity of $\hat{T}$.

For $n\geq q$ (where $q$ is the height of the tower), we can use
the same argument on the basis and get that
\begin{equation*}
  \norm{(\hat{T}^n \psi)_{|\Delta_0}}_{\boL_\tau}\leq B\left( \tau^n
  \norm{\psi}_{\boL_\tau}+\norm{\psi}_1\right).
\end{equation*}
If $x\in \Delta_k$, let $x'$ denote the corresponding point in the
basis of the tower. Then $\hat{T}^n
\psi(x)=\hat{T}^{n-k}\psi(x')$, and the inequality on the basis
gives
\begin{equation*}
  \norm{(\hat{T}^n \psi)_{|\Delta_k}}_{\boL_\tau}\leq B\left(
  \tau^{n-k}
  \norm{\psi}_{\boL_\tau}+\norm{\psi}_1\right).
\end{equation*}
As the tower is bounded, we finally get that
\begin{equation*}
  \norm{\hat{T}^n \psi}_{\boL_\tau}\leq \frac{B}{\tau^q}\left(
  \tau^n
  \norm{\psi}_{\boL_\tau}+\norm{\psi}_1\right).
\end{equation*}
This is a Doeblin-Fortet inequality (since the inclusion
$\boL_\tau \to L^1$ is compact). Thus, by Hennion's Theorem
(\cite{hennion}) the essential spectral radius of $\hat{T}$ is
bounded by $\tau$, and $1$ is an isolated eigenvalue of $\hat{T}$.
The constant functions are eigenfunctions, and ergodicity of $T$
shows that they are the only eigenfunctions (\cite[Theorem
1.4.8]{aaronson:book}). Moreover, 
there is no nilpotent part since $\|\hat{T}^n\|$ is
bounded.
\end{proof}

The following lemma will be useful later, applied to unbounded
functions in $L^2$. We recall that, if $d$ is an element of the
partition $\delta$, then $D_\tau h(d)$ is the least $\tau$-Hölder
constant of $h$ restricted to $d$.
\begin{lem}
\label{image_f_bornee} Assume that $T$ has the bounded tower big
image property. Let $h\in L^1(m)$ with $\sum_{d\in \delta}
m[d]D_\tau h(d)<\infty$. Assume that $h$ is bounded on the set of
points whose image is not in the basis of the tower, and that $h$
is uniformly $\tau$-Hölder on this set. Then $\hat{T}h \in
\boL_\tau$.
\end{lem}
\begin{proof}
For $\underline{d}=[d_0]$ of length $1$, the estimates given by
Lemma \ref{lemme_sur_Mdbarre} depend in fact only on
$D_\tau\psi(d_0)$. Thus, summing these estimates on the basis, and
using the big image property, we get
  \begin{equation*}
  \norm{(\hat{T}h)_{|\Delta_0}}_{\boL_\tau} \leq B\left(\sum_{d\in\delta}
  \tau m[d]D_\tau h(d)
  +\norm{h}_1\right)<\infty.
  \end{equation*}
If $x$ is in the floors, and $x'$ is just below $x$, then
$\hat{T}h(x)=h(x')$. As $h$ is bounded and Hölder on the set of
such points $x'$ by hypothesis, we get that $\hat{T}h$ is in
$\boL_\tau$.
\end{proof}

\subsection{Some technical estimates: the perturbative case}
\label{section_lemmes_perturbatifs} In this section,
$(Z,\boB,T,m,\delta)$ will satisfy the same assumptions as in the
previous section, i.e. $T$ is a Markov map, preserving the finite
measure $m$, with Hölder distortion. Moreover, $T$ has the bounded
tower big image property. Additionally, we will use a function
$h:Z\to \C$ which, on the set of points whose image is not in the
basis, will be bounded and uniformly $\tau$-Hölder, and satisfying
$\sum_{d\in\delta}m[d]D_\tau h(d)<\infty$.

Let $\underline{d}=[d_0]$ be a non-empty cylinder in $Z$ of length
$1$. For small $t$, we define a perturbation
$M_{\underline{d}}(t)$ of $M_{\underline{d}}$ by
  \begin{equation}
  \label{definit_Mt}
  M_{\underline{d}}(t)\psi(x)=
  \left\{\begin{array}{ll} e^{it h(d_0 x)} g_{m}(d_0 x) \psi(d_0
    x) & \text{ if $d_0x$ is defined}\\ 0 &\text{ otherwise}.
   \end{array}\right.
  \end{equation}

\begin{lem}
\label{lemme_perturbatif_sur_Mdbarre}
There exists a constant $B$ such that
$\forall \underline{d}=[d_0]$, $\forall t\in [-1,1]$,
  \begin{equation*}
  \norm{M_{\underline{d}}(t)-M_{\underline{d}}}_{\boL_\tau} \leq B
  \frac{m[\underline{d}]}{m[T \underline{d}]}
    \left( |t|D_\tau h(\underline{d})
    +\frac{1}{m[\underline{d}]}\int_{\underline{d}}
    |e^{it h}-1|\dd m \right).
  \end{equation*}
\end{lem}
\begin{proof}
Let $\psi\in \boL_\tau$. If $x,y$ are in the same element of the
partition, we have to study
$(M_{\underline{d}}(t)-M_{\underline{d}})\psi(x)
-(M_{\underline{d}}(t)-M_{\underline{d}})\psi(y)$. If $d_0 x$ is
not defined, then $d_0 y$ is not defined either, and there is
nothing to prove. Otherwise, denote $x'=d_0x$ and $y'=d_0y$. Then
  \begin{multline*}
  |(M_{\underline{d}}(t)-M_{\underline{d}})\psi(x)-
  (M_{\underline{d}}(t)-M_{\underline{d}})\psi(y)|
  \\
  \leq \left|e^{it h(x')}-1\right| |M_{\underline{d}}\psi(x)
  -M_{\underline{d}}\psi(y)|
  +|M_{\underline{d}}\psi(y)| \left|e^{it h(x')}-e^{it h(y')}\right|.
  \end{multline*}
We have
  \begin{equation*}
  \left|e^{it h(x')}-e^{it h(y')}\right|
  \leq |t| \left|h(x')-h(y')\right|
  \leq |t| D_\tau h(\underline{d}) \tau^{t(x,y)+1}.
  \end{equation*}
If $\eta=e^{it h}-1$, we see that $|\eta(x')|\leq |\eta(y')|+
|t|D_\tau h(\underline{d})$ on the cylinder $\underline{d}$,
whence $|e^{it h(x')}-1|\leq
\frac{1}{m[\underline{d}]}\int_{\underline{d}} |e^{it h}-1|\dd m
+|t| D_\tau h(\underline{d})$. The estimates on
$M_{\underline{d}}\psi$ from Lemma \ref{lemme_sur_Mdbarre} imply
that
  \begin{align*}
  |(M_{\underline{d}}(t)&-M_{\underline{d}})\psi(x)-
  (M_{\underline{d}}(t)-M_{\underline{d}})\psi(y)|
  \\
  \leq &
  \left(\frac{1}{m[\underline{d}]}\int_{\underline{d}} |e^{it h}-1|\dd m
  +|t|D_\tau h(\underline{d})\right) B \tau^{t(x,y)}
  \frac{m[\underline{d}]}{m[T\underline{d}]}\left(
  \tau D_\tau\psi+\frac{1}{m[\underline{d}]}
  \int_{[\underline{d}]} |\psi|
      \dd m\right)
  \\&\hphantom{=;}
  +
  B \frac{m[\underline{d}]}{m[T\underline{d}]}\left(
  \tau D_\tau\psi +\frac{1}{m[\underline{d}]} \int_{[\underline{d}]} |\psi|
      \dd m\right)
  |t|D_\tau h(\underline{d}) \tau^{t(x,y)}.
\end{align*}
This gives the estimate on $D_\tau \left(
(M_{\underline{d}}(t)-M_{\underline{d}})\psi\right)$, using the
fact that $\left(
  \tau D_\tau\psi+\frac{1}{m[\underline{d}]}
  \int_{[\underline{d}]} |\psi|
      \dd m\right) \leq 2 \norm{\psi}_{\boL_\tau}$. The bound on the
supremum norm is analogous.
\end{proof}

\begin{cor}
\label{controle_Rtp} Assume that $T$ has the bounded tower big
image property. Writing $R_t=\hat{T}(e^{ith}\cdot)$, then
$t\mapsto R_t$ is continuous at $0$, and
  \begin{equation*}
  \norm{R_t -R_0}_{\boL_\tau}=O(|t|+E(|e^{ith}-1|)).
  \end{equation*}
In particular, when $h\in L^1$, $\norm{R_t-R_0}=O(|t|)$.
\end{cor}
The notation $E$ denotes the expectation, i.e. the integral
 with respect to the measure
$\dd m$.

\begin{proof}
We will use that $R_t=\hat{T}(e^{ith}\cdot)=\sum
M_{\underline{d}}(t)$, where the sum is over the cylinders of length $1$.

To prove continuity of $t\mapsto R_t$ at $t=0$, we write
$R_t=1_{\Delta_0}R_t+(1-1_{\Delta_0})R_t$. On the basis, summing
the estimates of Lemma \ref{lemme_perturbatif_sur_Mdbarre} and
using the big image property, we get
 \begin{equation*}
  \norm{1_{\Delta_0}R_t-1_{\Delta_0}R_0}\leq B\left(|t|\sum m[\underline{d}]D_\tau
  h(\underline{d}) +\int |e^{it h}-1|\right)
  \end{equation*}
In the floors, if $x'$ is the point just below $x$, then
$R_t\psi(x)=e^{ith(x')}\psi(x')$, which yields that
$\norm{(1-1_{\Delta_0})R_t-(1-1_{\Delta_0})R_0}$ is bounded by the
$\boL_\tau$-norm of $e^{ith}-1$ restricted to the set of points
not returning to the basis at the next step. As $h$ is uniformly
bounded and Hölder on these points, we get that this
$\boL_\tau$-norm is $O(|t|)$. This proves that
$\norm{R_t-R_0}=O(|t|+\int |e^{ith}-1|)$. By the dominated
convergence theorem, the integral tends to $0$ when $t\to 0$, and
this proves the continuity of $R_t$.

Moreover, when $h$ is integrable, then $\int |e^{ith}-1|\leq
|t|\int |h|$.
\end{proof}

\subsection{The classical method: central limit theorem case}
\label{section_TCL_classique} In this section, we recall
classical results on the central limit theorem, obtained in the
uniform case by perturbing the transfer operator. For bounded $h$,
the method dates back to \cite{Nagaev}, and is developed in
\cite{rousseau-egele} and \cite{guivarch-hardy}. The result we
prove here is slightly more complicated to obtain, since the
function $h$ we will consider is generally unbounded, so that the
perturbation of the transfer operator is not analytic, and we can
not use an \emph{a priori} expansion of the eigenvalue. In view of
applications, we give the result for Markov maps with the bounded
tower big image property.

\begin{thm}
\label{TCL_classique} Let $(Z,\boB,T,m,\delta)$ be a Markov map,
preserving a finite ergodic measure $m$. Assume that $\log g_m$ is
$\tau$-Hölder for some $\tau<1$, and that $T$ has the bounded tower big
image property.

Let $h:Z\to \C$ with $\sum_{d\in\delta}m[d] D_\tau h(d)<\infty$
and $h\in L^2$, with $\int h \dd m=0$. Assume that, restricted
to the set of points that do not return to the basis of the tower
at the next iteration, $h$ is bounded and uniformly $\tau$-Hölder.

Denote by $R$ the transfer operator associated to $T$ (acting on
$\boL_\tau$), and $R_t \psi=R(e^{ith}\psi)$. Then, for small $t$,
$R_t$ has an eigenvalue $\lambda(t)$, close to $1$, satisfying
$\lambda(t)=1-\frac{\sigma^2}{2m(Z)}t^2+o(t^2)$. Write
$a=(I-R)^{-1}(Rh)$ (where we consider $I-R$ acting on the space of
functions of $\boL_\tau$ with $0$ average, on which it is
invertible) and $u=(I-R)^{-1}(h)$, then
  \begin{equation*}
  \sigma^2=\int h^2 \dd m+2\int ah\dd m
  =\int R(u^2)-(Ru)^2 \dd m.
  \end{equation*}

Moreover, $\sigma^2=0$ if and only if there exists a function
$\psi\in \boL_\tau$ with $h=\psi\circ T-\psi$.
\end{thm}
Note that Lemma \ref{image_f_bornee} ensures that $Rh\in
\boL_\tau$, which implies that it makes sense to consider
$a=(I-R)^{-1}(Rh)$, since $I-R$ is invertible on the functions of
$\boL_\tau$ of zero integral, according to Corollary
\ref{big_donc_vpisolee}. However, $u=(I-R)^{-1}(h)$ is \emph{a
priori} not defined, since $h\not\in \boL_\tau$. In fact, we
simply set $u=h+a$, and indeed $(I-R)u=h-Rh+ (I-R)a=h$.

\begin{proof}
We can replace $m$ by $m/m(Z)$ and assume that $m$ is a
probability measure (the greater generality in the statement of the
theorem will be useful in the
applications).

Corollary \ref{big_donc_vpisolee} gives that $R=\hat{T}$ is
continuous on $\boL_\tau$, and $1$ is a simple isolated
eigenvalue. According to Corollary \ref{controle_Rtp}, $t\mapsto
R_t$ is continuous at $0$. As simple isolated eigenvalues depend
continuously on the operator, $R_t$ has a unique eigenvalue
$\lambda(t)$ close to $1$ for $t$ small. Write $P_t$ for the
corresponding spectral projection, and $\xi_t$ for the
eigenfunction (with $\int \xi_t=1$).

Since eigenvalues and  eigenfunctions depend
holomorphically on operators (Lipschitz would be enough),
$\norm{\xi_t-1}=O(\norm{R_t - R})=O(|t|)$ using Corollary
\ref{controle_Rtp}. In the same way, $\lambda(t)=1+O(t)$.

Note for future use that
  \begin{equation}
  \label{fct_caract}
  E(e^{ith})=1-\frac{t^2}{2}\int h^2 \dd m +o(t^2).
  \end{equation}
Indeed, if we consider $\Phi(t)=E(e^{ith})$ as the characteristic
function of the random variable $h$, it is classical that if
$h\in L^2$, then $\Phi$ is $C^2$ with $\Phi(0)=1$,
$\Phi'(0)=E(h)$ and $\Phi''(0)=E(h^2)-E(h)^2$ (\cite[Corollary to
Lemma XV.4.2]{feller:2}).

As $\lambda(t) \xi_t=R_t \xi_t$, we get after integration that
  \begin{equation}
  \begin{split}
  \label{estime_lambdat_direct}
  \lambda(t)&=\int R_t \xi_t
  =\int (R_t-R) (\xi_t-1)+\int (R_t-R)(1)+\int R \xi_t
  \\
  &=\int (R_t-R)(\xi_t-1)+E(e^{ith}-1)+1.
  \end{split}
  \end{equation}
As $R_t-R=o(1)$ and $\xi_t-1=O(t)$, we get that $\lambda(t)
=1+o(t)$ (using Equation \eqref{fct_caract} to estimate
$E(e^{ith}-1)$).

Thus,
  \begin{align*}
  \frac{\xi_t-\xi_0}{t}&
  =\frac{\lambda(t) \xi_t -\xi_0}{t}+o(1)
  =\frac{R_t \xi_t -R_0 \xi_0}{t}+o(1)
  \\
  &=(R_t-R_0)\frac{\xi_t-\xi_0}{t}+R_0 \frac{\xi_t-\xi_0}{t}
  +\frac{R_t-R_0}{t} \xi_0 +o(1).
  \end{align*}
But $(R_t-R_0)(\xi_0)=R(e^{ith}-1)$. Moreover,
$\frac{\xi_t-\xi_0}{t}$ is bounded and $R_t-R_0$ tends to $0$,
whence $(R_t-R_0)\frac{\xi_t-\xi_0}{t}=o(1)$. Thus,
  \begin{equation*}
  (I-R) \frac{\xi_t-\xi_0}{t}=R\left(\frac{e^{ith}-1}{t}\right)+o(1).
  \end{equation*}
The sequence $(\xi_t-\xi_0)/t$ is bounded in $\boL_\tau$. By
compactness, we can extract a subsequence converging in $L^2$ to a
function $ia$, which will still be in $\boL_\tau$ (we have
multiplied by $i$ for convenience). In $L^2$,
$\frac{e^{ith}-1}{t}\to ih$ by the dominated convergence theorem
(as $\left|\frac{e^{ith}-1}{t}\right|^2\leq |h|^2$ integrable),
which implies that $(I-R)a=R(h)$. Moreover, according to Lemma
\ref{image_f_bornee}, $Rh \in \boL_\tau$. We also have $\int
Rh=\int h=0$. As $I-R$ is invertible on the set of functions in
$\boL_\tau$ with zero integral, we get $a=(I-R)^{-1}(Rh)$. Thus,
$(\xi_t-\xi_0)/t$ has a unique cluster value, whence it converges
(in $L^2$) to $ia=i(I-R)^{-1}(Rh)$.

By continuity of the product $L^2 × L^2 \to L^1$, we obtain the
convergence of $\int \frac{e^{ith}-1}{t} \frac{\xi_t-\xi_0}{t}$ 
to $\int (ih)\cdot (ia)$. Thus,
  \begin{equation*}
  \int (R_t-R)(\xi_t-1)
  =\int R( (e^{ith}-1) (\xi_t-1))
  =\int (e^{ith}-1)(\xi_t-1)
  =-t^2 \int ha +o(t^2).
  \end{equation*}
Equation \eqref{estime_lambdat_direct} then gives that
  \begin{align*}
  \lambda(t)=1+E(e^{ith}-1)+\int (R_t-R)(\xi_t-1)
  =1-\frac{t^2}{2}\int h^2 -t^2 \int ha +o(t^2).
  \end{align*}
We transform this result a little to obtain the other claimed
expression for $\lambda(t)$. Recall that $h=u-Ru$ and $a=u-h$.
  \begin{align*}
  -\lambda''_{|t=0}&=\int h^2+2\int a h
  =\int h^2 +2(u-h)h
  =\int (u-Ru)^2+\int 2Ru\cdot (u-Ru)
  \\&
  =\int (u+Ru)(u-Ru)
  =\int u^2-\int (Ru)^2
  =\int R(u^2)-(Ru)^2.
  \end{align*}
We have $R(u^2)(x)=\sum_{Ty=x}g_m(y)u^2(y)$ and
$(Ru)^2(x)=\left(\sum_{Ty=x} g_m(y)u(y)\right)^2$. Moreover,
$\sum_{Ty=x}g_m(y)=1$ since $R1=1$. The convexity of $w\mapsto
w^2$ then gives that $(Ru)^2(x)\leq R(u^2)(x)$, which implies that
$-\lambda''\geq 0$. We can thus write $\lambda''=-\sigma^2$.

Finally, if $\sigma^2=0$, then there is equality in the convexity
inequality, whence all the $u(y)$ are equal (since $g_m$ is everywhere
nonzero as $
\log g_m$ is defined everywhere). Writing $\psi=Ru$ (with $\psi\in\boL_\tau$
since $Ru=Ra+Rh$) we get $u=\psi\circ T$.
As $h=u-Ru$, this gives $h=\psi\circ T-\psi$.
\end{proof}

\subsection{The classical method: stable laws case}
\label{section_stable_classique}

The following theorem has essentially been proved by Aaronson and
Denker in \cite{aaronson_denker}.
\begin{thm}
\label{loi_stable_classique} Let $(Z,\boB,T,m,\delta)$ be a Markov
map preserving the finite ergodic measure $m$. Assume that
$\log g_m$ is $\tau$-Hölder for some $\tau<1$ and that $T$ has the
bounded tower big image property.

Let $p\in (0,1)\cup (1,2)$ and $h$ be a function on $Z$ with
$\sum_{d\in\delta}m[d] D_\tau h(d)<\infty$, such that, on the
elements of the partition not returning to the basis of the tower
at the next iteration, $h$ is bounded and uniformly Hölder.
Moreover, assume that $m[h>x]=(c_1+o(1))x^{-p}L(x)$ and
$m[h<-x]=(c_2+o(1))x^{-p}L(x)$ for constants $c_1,c_2\geq 0$ with
$c_1+c_2>0$, and a slowly varying function $L$. Finally, assume
$\int h=0$ if $p>1$.

Let $R$ be the transfer operator associated to $T$, acting on
$\boL_\tau$, and $R_t \psi=R(e^{ith}\psi)$. Then, for small enough
$t$, $R_t$ has a unique eigenvalue close to $1$, and this
eigenvalue satisfies
  \begin{equation*}
  \lambda(t)=1-\frac{c}{m(Z)}|t|^p L(1/|t|)+i\frac{c}{m(Z)}|t|^p L(1/|t|)\beta
  \sgn(t)\tan\left(\frac{p\pi}{2}\right)+o(|t|^pL(1/|t|))
  \end{equation*}
for constants $\beta=\frac{c_1-c_2}{c_1+c_2}$ and
$c=(c_1+c_2)\Gamma(1-p)\cos\left(\frac{p\pi}{2}\right)$.
\end{thm}
\begin{proof}
Replacing $m$ by $m/m(Z)$, we can assume that $m$ is a probability measure.

This theorem is essentially Theorem 5.1 in \cite{aaronson_denker},
with some differences in the hypotheses: we use the bounded tower
big image property, whence their proof has to be slightly
modified, as we explain now.

Corollary \ref{big_donc_vpisolee} gives that $1$ is a simple
eigenvalue of $R$ acting continuously on $\boL_\tau$. Corollary
\ref{controle_Rtp} ensures that $t\mapsto R_t$ is continuous at
$0$, which implies that $R_t$ also has a simple isolated
eigenvalue $\lambda(t)$ close to $1$. The proof of Theorem 5.1 in
\cite{aaronson_denker} applies then literally if we can prove
that, for $p>1$, the eigenfunction $\xi_t$ of $R_t$ (normalized to
have integral $1$) satisfies $\norm{\xi_t-1}=O(|t|)$.

As $\xi_t$ is an eigenfunction of $R_t$ and eigenfunctions depend
holomorphically on operators, $\norm{\xi_t-1}=O(\norm{R_t-R})$.
Corollary \ref{controle_Rtp} gives that this is a $O(|t|)$, which
concludes the proof.
\end{proof}

An analogous theorem for $p=2$ and $h\not\in L^2$ may be proved
using \cite{aaronson_denker:central}.

\section{Proof of the main theorems}
\label{section:proof_thm}

In this section, we will prove Theorems \ref{TCL_markov} and
\ref{thm_loi_stable_markov} using the results of Sections
\ref{section:abstract_result} and \ref{section:uniform_case}.
$(X,\boB,T,m,\alpha)$ will be a probability preserving
topologically mixing Markov map, such that the induced map on
$Y\subset X$ has the big image property and Hölder distortion. The
function $\phi=\phi_Y$ will denote the first return time from $Y$
to itself.

Note first that, since $T$ is topologically mixing, $T_Y$ is topologically
transitive. Theorem 4.6.3 in \cite{aaronson:book} implies that $T_Y$ is
ergodic, whence $T$ is also ergodic.

\subsection{Construction of an extension}
\label{construit_Kakutani}
In order to change variables
between $\{x \tq \phi(x)>i\}$ and its image by $T^i$, we have to
construct an extension of the system since it is possible that two
different points in $Y$ are sent on the same image in $X$ before
they return to $Y$.

We use a tower extension, sometimes called a Kakutani tower (the
towers in Section \ref{section:uniform_case} will actually correspond to
making a finite height cutoff of the Kakutani tower). It is built
as follows: we set $X'=\{(x,i)\tq x\in Y, 0\leq i<\phi(x)\}$. $Y$
is identified with $Y':=\{(x,0) \tq x\in Y\}\subset X'$. The space
$X'$ is endowed with a measure $m'$, equal to $m$ on $Y'$ and then
lifted in the tower, i.e. if $A\subset Y$ and $\phi>i$ on $A$,
then $m'(A\times\{i\})=m(A)$.

We define $T'$ on $X'$ by $T'(x,i)=(x,i+1)$ if $i<\phi(x)-1$, and
$(T^{\phi(x)}x,0)$ otherwise. We have a projection $\pi:X'\to X$
defined by $\pi(x,i)=T^i x$ (which is not necessarily injective). The
following diagram is then commutative, and all applications are
measure preserving:
  \begin{equation*}
  \xymatrix{
  X' \ar[d]_\pi \ar[r]^{T'} &
       X' \ar[d]^\pi \\
  X \ar[r]^{T} & X}
  \end{equation*}

All this is classical material, but I have not been able to
locate a satisfying reference in the literature. Thus, we give for
completeness a proof of the fact that the map $\pi$ is
measure-preserving.
\begin{prop}
$\forall K\subset X$, we have $m'(\pi^{-1} K)=m(K)$.
\end{prop}
\begin{proof}
Let $A_0=K- Y$, $B_0=K\cap Y$. We define inductively
$A_k=T^{-1}(A_{k-1})- Y$ and $B_k=T^{-1}(A_{k-1})\cap Y$: the points of
$A_k$ have not been treated up to time $k$, while the points in $B_k$
are treated at time $k$. Then we have $\pi^{-1}(K)=\{ (x,i) \tq i\in
\N, x\in
B_i\}$. Thus, $m'(\pi^{-1}K)=\sum m(B_k)$.

As $m(K)=m(A_0)+m(B_0)$ and, by construction,
$m(A_{k-1})=m(T^{-1}A_{k-1})=m(A_k)+m(B_k)$, we get that for any $n$,
$m(K)=m(B_0)+\ldots +m(B_n)+m(A_n)$. To conclude, we just have to show
that $m(A_n)\to 0$ when $n\to \infty$.

Let $C_n=\{x \tq \forall i\leq n, T^i x\not \in Y\}$. Then
$A_n\subset C_n$. Moreover, $C_n$ is decreasing. Let $C=\bigcap
C_n$. It suffices to show that $m(C)=0$, since this will prove
that $m(C_n)\to 0$. As $C \subset T^{-1}C$ and $T$ is measure
preserving, $C=T^{-1}C$ mod $0$. By ergodicity of $T$, $m(C)=0$ or
$1$. As $C\cap Y=\emptyset$ and $m(Y)>0$, we obtain $m(C)=0$.
\end{proof}

If $\chi$ is a function on $X$, we define $\chi'$ on $X'$ by
$\chi'=\chi \circ \pi$. Then $\int_{X'} \chi' \dd m'=\int_X
\chi \dd m$.

To construct a partition on $X'$, we consider first the partition
$\delta$ on $Y$ for which the induced map $T_Y$ is Markov. Then,
we lift it in the tower, i.e. the partition at height $i$ is the
same as the partition of the projection of this floor on the
basis. Then $T'$ is Markov for this new partition $\alpha'$. Note
that the partition for the induced map $T'_{Y'}$ is then simply
the restriction of the partition $\alpha'$ to the basis of the
tower.

Moreover, the inverse of the jacobian of $T'_{Y'}$ is the same as
the inverse of the jacobian of $T_Y$, whence it is locally
$\theta$-Hölder (note that the separation times are the same in
$Y$ and in $Y'$).

Finally, $T'$ is still topologically mixing: if $a\in \delta$,
then one of its iterates $T^n a$ contains an element $b$ of
$\alpha$. If $c\in \delta$, it is contained in some $c'\in
\alpha$, and there exists $N$ such that $\forall p\geq N$,
$c'\subset T^p(T^n a)$. This proves that $\forall p\geq n+N$, $c
\subset T^p(a)$. This easily implies topological mixing for $T'$.

In fact, all the hypotheses on
$(X,T,f,f_Y)$ are carried over to $(X',T',f',f_Y')$. If we prove the
central limit theorem or the convergence to a stable law on $X'$, this
will give the same result on $X$ since $(S_n f)'=S_n(f')$. Thus, it is
sufficient to prove the result for $X'$.

From this point on, we will thus work in a tower. It is Markov,
and the returns to the basis have a big image (but are not
necessarily surjective). We will first prove Theorem
\ref{TCL_markov} on the Central Limit Theorem, and we will
indicate in a last section the modifications to be done for
Theorem \ref{thm_loi_stable_markov} (the stable case).

\subsection{Local result}

The hypotheses of the main theorems are tailor-made so that
induction on the basis $Y$ of the tower has good properties.
However, we will have to induce on larger parts of the space. In
this section, we describe properties of the induced maps on
bounded parts of the tower. In this section, the hypotheses will
be those of Theorem \ref{TCL_markov}.

\subsubsection{The space}
Let $q\geq 1$. We write $Z=Z_q$ for the union of the first $q$
floors of the Kakutani tower from Section
\ref{construit_Kakutani}. We will induce on $Z$. We will write
$\delta$ for the partition on $Z$ for which the induced map $T_Z$
is Markov -- this is in fact simply the restriction of the
original partition $\alpha$ to $Z$, because of the particular
combinatorics of the tower. For $k<q$, we will also denote by
$\Delta_k$ the $k^{\text{th}}$ floor of the tower. Note that the
hypotheses guarantee that the function $f$ is bounded and Hölder
on $Z$.

If $s$ is the separation time with respect to the basis $Y$ of the
tower, and $t$ is the separation time for $T_Z$, then an iteration
of $T_Y$ corresponds to at most $q$ iterations of $T_Z$, which
implies that $t(x,y)\leq q s(x,y)$. Thus, if $\psi$ satisfies
$|\psi(x)-\psi(y)|\leq C\theta^{s(x,y)}$, it will
  satisfy $|\psi(x)-\psi(y)|\leq C\tau^{t(x,y)}$ for $\tau=\theta^{1/q}$.

We write $g_{m_Z}=\frac{\dd m_Z}{\dd m_Z\circ T_Z}$. For $(x,i)\in
Z$, if $i<\phi(x)-1$ then $g_{m_Z}(x)=1$, and if $i=\phi(x)-1$
then $g_{m_Z}(x,i)=g_{m_Y}(x,0)$ (the inverse of the jacobian of
$T_Y$ with respect to $m_Y$). As $\log g_{m_Y}$ is locally
$\theta$-Hölder (for the separation time $s$), we obtain that
$\log g_{m_Z}$ is locally $\tau$-Hölder (for the separation time
$t$).

Since $T$ is ergodic, so is $T_Z$ (\cite[Proposition 1.5.2]{aaronson:book}).
Moreover, $T_Z$ has the bounded tower big image property. Thus,
the transfer operator $\hat{T_Z}$ acts continuously on
$\boL_\tau(Z)$ according to Corollary \ref{big_donc_vpisolee}, and
has a simple isolated eigenvalue at $1$.

\subsubsection{First return transfer operators}

Following \cite{sarig:decay}, we define operators $R_n$ and $T_n$:
$R_n$ is a first return transfer operator, i.e. it sees only the first
returns to $Z$ at time $n$, while $T_n$ sees all returns at time $n$.
They are defined by
$R_n \psi= 1_Z \hat{T}^n(1_{\{\phi_Z=n\}}
\psi)$ (where $\phi_Z$ is the first return time from $Z$ to itself)
and $T_n \psi= 1_Z \hat{T}^n (1_Z \psi)$. They can also be
written in the following way: if $x\in Z$, then $R_n\psi(x)=\sum
g_m^{(n)}(y)\psi(y)$ where the sum is over all points $y$ in $Z$ such that
$Ty,\ldots, T^{n-1}y \not \in Z$ and $T^n y=x$, while $T_n
\psi(x)=\sum g_m^{(n)}(y)\psi(y)$ where the sum is over all points
$y\in Z$ with $T^n y=x$.

If $y\in Z$ with $T^n y=x$, we can consider all its iterates between
$0$ and $n$ that fall into $Z$. To go from one of these iterates to
the following corresponds to iterating one of the $R_k$s.
Thus, we get
$T_n=\sum_{k_1+\ldots+k_l=n} R_{k_1}\ldots R_{k_l}$: this is the
\emph{renewal equation}, which will make it possible to understand the
$T_n$ if the $R_n$ are well understood.

\begin{lem}
\label{lemme_Rn_continu}
The operators $R_n$ and $T_n$ act
continuously on $\boL_\tau(Z)$, and $\norm{R_n}_{\boL_\tau(Z)}
=O(m[\phi_Z=n])$ (where $\phi_Z$ is the return time from $Z$ to
itself).
\end{lem}
\begin{proof}
As $\hat{T_Z}$ acts continuously on $\boL_\tau$ and
$R_1\psi=\hat{T_Z}(1_{\{\phi_Z=1\}} \psi)$, the operator $R_1$ is
  continuous.

For $n\geq 2$, we estimate $\norm{R_n}$ using the notations of
Section \ref{section_quelques_lemmes} (applied to $Z$ and $T_Z$):
  \begin{equation*}
  R_n=\sum_{\underline{d}=[d_0]_Z,
  d_0=[a_0,\ldots,a_{n-1},Z]} M_{\underline{d}}.
  \end{equation*}
For $n\geq 2$, the cylinders appearing in $R_n$ correspond to
points coming back to $Z$ in $n$ steps: they get out of the tower
$Z$ through the top, spend $n-1$ iterates in the higher floors of the
Kakutani tower, and then get back to the basis of the tower, with
a big image. Thus, we can forget the terms $m[T_Z d_0]$ for these
cylinders in the estimates of Lemma \ref{lemme_sur_Mdbarre}, and
we get $\norm{R_n}\leq B(1+\tau) \sum m[d_0]=B(1+\tau)
m[\phi_Z=n]$.

The $R_n$ are continuous and $T_n=\sum_{k_1+\ldots+k_l=n} R_{k_1}\ldots
R_{k_l}$, whence $T_n$ is also continuous.
\end{proof}

\begin{lem}
\label{lemme_inversible_hors_de_1}
For $z\in \overline{\D}$, write
$R(z)=\sum R_n z^n$. Then $R(1)$ acting on $\boL_\tau(Z)$ has a
simple isolated eigenvalue equal to $1$, and the corresponding
spectral projection is given by $\psi \mapsto \frac{1}{m[Z]}
\int_Z \psi \dd m$.

Moreover, for $z\not=1$, $I-R(z)$ is invertible.
\end{lem}
Note that, since $\norm{R_n}=O(m[\phi_Z=n])$ is summable by Kac's
Formula, the series $R(z)$ is converging for all $z\in \overline{\D}$.
\begin{proof}
Let $z\in\overline{\D}$. We will show a Doeblin-Fortet inequality
for $R(z)^s$ when $s\geq q$. Let $\psi\in \boL_\tau(Z)$. Then
  \begin{equation*}
  R(z)^s=\sum_{\underline{d}=[d_0,\ldots,d_{s-1}]_Z} z^{l(\underline{d})}
  M_{\underline{d}}
  \end{equation*}
where $l(\underline{d})$ is the sum of the lengths of the $d_i$s
seen as cylinders in $X$. In particular, $l(\underline{d})\geq s$,
whence $|z|^{l(\underline{d})}\leq |z|^s$. On the basis, summing
the inequalities given by Lemma \ref{lemme_sur_Mdbarre} (and using
the big image property for the returns to the basis), we obtain
that $\norm{(R(z)^s\psi)_{|\Delta_0}}_{\boL_\tau}
 \leq B |z|^s\left( \tau^s \norm{\psi}_{\boL_\tau}+
  \int |\psi|\right)$.

On $\Delta_k$ the set of points at height $k$, $R(z)^s \psi(x)=z^k
R(z)^{s-k} \psi(x')$ where $x'$ is the point corresponding to $x$
but in the basis. Thus, $\norm{\left(R(z)^s
  \psi\right)_{|\Delta_k} } \leq B |z|^s \left( \tau^{s-k}
\norm{\psi}_{\boL_\tau}+
  \int |\psi|\right)$.

Summing on all the floors, we have proved that
  \begin{equation*}
  \forall s\geq q, \forall \psi\in \boL_\tau,
  \norm{R(z)^s \psi}_{\boL_\tau}
  \leq \frac{B}{\tau^q} |z|^s \left( \tau^s \norm{\psi}_{\boL_\tau}+
    \int |\psi|\right).
  \end{equation*}

As the injection from $\boL_\tau$ to $L^1$ is compact, this is a
Doeblin-Fortet inequality, which implies that the essential spectral
radius of $R(z)$ is $\leq \tau |z|<|z|$. Moreover, this
inequality also gives that the spectral radius of $R(z)$ is $\leq |z|$.

If $|z|<1$, then $R(z)$ has spectral radius $\leq |z|<1$, whence
$I-R(z)$ is invertible.

If $z=1$, as $R(1)$ counts the first returns to $Z$, it is not hard to check
that $R(1)=\hat{T_Z}$ is the transfer operator associated to
$T_Z$. Hence, Corollary \ref{big_donc_vpisolee}
ensures that it has a simple isolated
eigenvalue at $1$.

If $|z|=1$ but $z\not=1$, we have to show that $I-R(z)$ is
invertible. We write $z=e^{it}$ for some $0<t<2\pi$.
Since the essential spectral radius of $R(z)$ is $<1$, it
is enough to show that $1$ is not an eigenvalue of $R(z)$. Suppose on
the contrary that $R(z)a=a$ for some nonzero $a\in \boL_\tau(Z)$.

For $\psi,\xi \in L^2(m_Z)$, write $\langle \psi,\xi\rangle=\int
\overline{\psi}\xi \dd m_Z$. Define the operator
$W:L^\infty(m_Z)\to L^\infty(m_Z)$ by
$W\psi=e^{-it\phi_Z}\psi\circ T_Z$. As
$R(z)\xi=R(1)(e^{it\phi_Z}\xi)$, the operator $W$ satisfies
    \begin{equation*}
    \langle \psi,R(z)\xi\rangle
    =\int \overline{\psi}\,R(z)\xi
    =\int \overline{\psi}\,R(1)(e^{it\phi_Z}\xi)
    =\int \overline{\psi}\circ T_Z\, e^{it\phi_Z}\xi
    =\int \overline{W\psi}\cdot \xi
    =\langle W\psi,\xi\rangle.
    \end{equation*}
We show that $a$ is an eigenfunction of $W$ for the eigenvalue $1$:
    \begin{align*}
    \norm{Wa-a}_2^2&
    =\norm{Wa}_2^2-2\Ree\langle Wa,a\rangle+\norm{a}_2^2
    =\norm{Wa}_2^2-2\Ree\langle a,R(z)a\rangle +\norm{a}_2^2
    \\&
    =\norm{Wa}_2^2-2\Ree\langle a,a\rangle+\norm{a}_2^2
    =\norm{Wa}_2^2-\norm{a}_2^2.
    \end{align*}
As $T_Z$ preserves the measure $m_Z$, we have
$\norm{Wa}_2^2=\int |a|^2\circ T_Z=\int |a|^2=\norm{a}_2^2$,
which gives $\norm{Wa-a}_2^2=0$. Hence, the function $Wa-a$ is
zero $m_Z$-almost everywhere. As $a\in \boL_\tau$ and $m_Z$
is nonzero on every cylinder, the function $a$ is continuous, whence
$Wa-a=0$ everywhere, i.e. $e^{-it \phi_Z} a\circ T_Z=a$.

Let $b$ be the restriction of $a$ to the basis $Y$ of the tower.
Then the previous equality implies that $e^{-it \phi_Y} b\circ
T_Y=b$. Taking the modulus, ergodicity of $T_Y$ gives that $|b|$
is constant almost everywhere, hence everywhere by continuity. As
$b\not\equiv 0$, this constant is nonzero, and we get
$e^{-it\phi_Y}=b/b\circ T_Z$. We can apply Theorem 3.1. in
\cite{aaronson_denker} since $T_Y$ has the big image property
(this is not the case for $T_Z$, which is why we have to restrict
$a$ to the basis of the tower). Writing $\beta$ for the induced
partition on the basis, this theorem gives that $b$ is
$\beta^*$-measurable, where $\beta^*$ is the smallest partition
such that $\forall d\in \beta$, $T_Y d$ is contained in an atom of
$\beta^*$. In particular, $b$ is constant on each set of $\beta$.

Let $d\in \beta$. On $[d]$, $b$ is equal to a constant $c$. As
$T$ is topologically mixing, there exists $N$ such that, $\forall
n\geq N$, $[d]\subset T^n[d]$. Let $n\geq N$, and $x\in [d]$ be
such that $T^n x\in [d]$. Let $T^{k_1}x, T^{k_2}x,\ldots,T^{k_p}x$
be the successive returns of $x$ to $Y$, with $k_p=n$. Then
$T^nx=T_Y^px$ and $n=\sum_{k=0}^{p-1} \phi_Y(T_Y^k
x)$. Thus,
  \begin{equation*}
  e^{-itn}=e^{-it \sum_{k=0}^{p-1} \phi_Y(T_Y^k x)}
  =\frac{b(x)}{b(T_Y x)} \frac{b(T_Y x)}{b(T_Y^2 x)} \cdots
  \frac{b(T_Y^{p-1}x)}{b(T_Y^px)}
  =\frac{b(x)}{b(T^n x)}=\frac{c}{c}=1.
  \end{equation*}
This is true for any $n\geq N$. Taking for example $n=N$ and $N+1$ and
quotienting, we obtain $e^{it}=1$, which is a contradiction and
concludes the proof.
\end{proof}

\begin{lem}
\label{lemme_valeur_de_mu} Writing $R'(1)=\sum_{n=1}^\infty nR_n$
and $P$ for the spectral projection associated to the eigenvalue
$1$ of $R(1)$, we have $PR'(1)P=\mu P$ for $\mu=\frac{1}{m[Z]}>0$.
\end{lem}
\begin{proof}
Since $Pf=\frac{1}{m[Z]}\int_Z f\dd m$, we have
$PR_nP=\frac{m[\phi_Z=n]}{m[Z]}P$, whence
  \begin{equation*}
  PR'(1)P=\left(\sum_{n=1}^{\infty} n m[\phi_Z=n]\right) \frac{1}{m[Z]}P
  =\frac{1}{m[Z]}P
  \end{equation*}
by Kac's Formula.
\end{proof}

\subsubsection{Perturbation of the transfer operators}
Recall that we are in the setting of Theorem \ref{TCL_markov}: $f$
is a function on $X$ such that the induced function
$f_Y=\sum_{0}^{\phi_Y-1} f\circ T^k$ is H\"older and square
integrable on $Y$.

Recall also that $\phi_Z$ is the first return time to our truncated
tower $Z$, and define a function $f_Z$ on $Z$ by
$f_Z(x)=\sum_{k=0}^{\phi_Z(x)-1} f(T^kx)$. If $x$ does not go out
of the bounded tower $Z$ through the top, we have $f_Z(x)=f(x)$,
while otherwise $f_Z(x)$ is the sum of the $f(y)$ for $y$ above
$x$ in the tower. Since $f\in L^1(m)$, we also have $f_Z\in
L^1(m_Z)$ and $\int_Z f_Z=\int_X f=0$. This function $f_Z$ is
interesting since we can study the Birkhoff sums of $f$ for $T$ by
looking at the Birkhoff sums of $f_Z$ for $T_Z$: if $x\in Z$ and
$T^n x\in Z$, and if $k$ is the number of returns of $x$ to $Z$
between times $0$ and $n$, then $\sum_{i=0}^{n-1} f(T^i
x)=\sum_{j=0}^{k-1}f_Z(T_Z^j x)$.

We show that $f_Z\in L^2(Z)$ (recall that we are proving the
Central Limit Theorem, i.e. we assume that $f_Y\in L^2(Y)$). Let
$\Delta=\Delta_{q-1}$ be the highest floor of $Z$, and $\Lambda$
its projection on the basis. Then $f_Z$ is bounded on $Z-\Delta$.
Moreover, if $x\in \Delta$ and $x'$ is the corresponding point in
$\Lambda$, then $f_Y(x')=f_Z(x)+\sum_0^{q-2} f(T^j x')$, and the
sum is bounded by $q\norm{f_{|Z}}_\infty$, whence
$f_Y(x')=f_Z(x)\pm C$. As ${f_Y}_{|\Lambda}$ is in $L^2$, we
obtain also ${f_Z}_{|\Delta}\in L^2$.

We now perturb the first return transfer operators, setting
$R_n(t)=R_n(e^{itf_Z} \cdot)$ and $T_{n,t}=R_n(e^{itS_nf}\cdot)$.
Since the Birkhoff sums of $f$ for $T$ and of $f_Z$ for $T_Z$ are
equal, the renewal equation still holds for these operators, i.e.
$T_{n,t}=\sum_{k_1+\ldots+k_l=n} R_{k_1}(t)\ldots R_{k_l}(t)$.
This equation can also be written as $\sum T_{n,t}z^n=\left(I-\sum
R_n(t) z^n\right)^{-1}$.

\begin{lem}
\label{lemme_rn_sommable}
There exist constants $r_n$ and
$\epsilon_0>0$ such that, $\forall n\geq 1, \forall t \in
[-\epsilon_0,\epsilon_0]$, $\norm{R_n(t)}\leq r_n$, and satisfying
$\sum n r_n<\infty$.

Moreover, for every $n\geq 1$, the map $t\mapsto R_n(t)$ is
continuous at $t=0$.
\end{lem}
\begin{proof}
Let $n\geq 2$. We have $R_n(t)=\sum M_{\underline{d}}(t)$ where
the sum is over all $\underline{d}=[d_0]_Z$ with
$d_0=[a_0,\ldots,a_{n-1},Z]$ (we use here the operators
$M_{\underline{d}}(t)$ introduced in \eqref{definit_Mt}, for the
function $h=f_Z$). The returns to $Z$ do not take place
at the first iteration, so they have to be returns to the basis,
whence they have a big image. Summing the estimates given by Lemma
\ref{lemme_perturbatif_sur_Mdbarre}, we obtain
  \begin{align*}
  \norm{R_n(t)}_{\boL_\tau}&\leq\norm{R_n}_{\boL_\tau}+ B
  \left(\sum_{d\in\delta,\phi_Z(d)=n}m[d]D_\tau f_Z(d)
  +\int_{\{\phi_Z=n\}} |e^{itf_Z}-1|\dd m \right)
  \\&
  \leq C m[\phi_Z=n]+B\sum_{d\in\delta,\phi_Z(d)=n}m[d]D_\tau f_Z(d)
  + 2B m[\phi_Z=n]=:r_n.
  \end{align*}
We know that $\sum n m[\phi_Z=n]=m(X)<\infty$, by Kac's formula.
Thus, to show that $\sum n r_n<\infty$, it is sufficient to prove
that $\sum_{d\in\delta}\phi_Z[d]m[d]D_\tau f_Z(d)<\infty$. Note
that this inequality holds for $f_Y$ by assumption.

Let $A$ be the set of elements of the partition of $Z$ that do not
return to the basis of $Z$ at the next iteration, and $B$ the
set of those that come back to the basis. If $d\in A$, then
$\phi_Z=1$ on $d$, and $f_Z=f$ is uniformly $\tau$-H\"older,
whence
  \begin{equation*}
  \sum_{d\in A} \phi_Z[d]m[d] D_\tau f_Z(d) \leq C\sum_{d\in A}
  m[d]\leq Cm[Z] <\infty.
  \end{equation*}
Let $d\in B$ and $x,y\in d$. Let $x'$ and $y'$ be the points in
the basis corresponding to $x$ and $y$, and $d'$ be the
corresponding partition element in the basis: $x=T^i x'$ and
$y=T^i y'$ for some $i<q$. Then $f_Y(x')=\sum_0^{i-1} f(T^j
x')+f_Z(x)$, and we have the same expression for $y$, whence
  \begin{align*}
  |f_Z(x)-f_Z(y)|&\leq |f_Y(x')-f_Y(y')|+\sum_{j=0}^{i-1} |f(T^jx')-f(T^jy')|\\
  &\leq C D_\theta f_Y(d') \theta^{s(x',y')}+ D q \theta^{s(x',y')}
  \leq (CD_\theta f_Y(d')+Dq)\tau^{t(x,y)}
  \end{align*}
using $s(x',y')=s(x,y)\geq t(x,y)/q$ and $\tau=\theta^{1/q}$.
Since $\phi_Z(d)\leq \phi_Y(d')$, we get
  \begin{equation*}
  \sum_{d\in B} \phi_Z[d]m[d]D_\tau f_Z(d)
  \leq \sum_{d'} \phi_Y[d']m[d'](CD_\theta
  f_Y(d')+Dq)<\infty
  \end{equation*}
by Kac's formula and the assumption on $f_Y$. This proves that
$\sum nr_n<\infty$.

For $n=1$, $R_1(t)=R_t(1_{\{\phi_Z=1\}}\cdot)$ and $R_t$
depends continuously on $t$ at $t=0$ according to Corollary
\ref{controle_Rtp}.

Finally, to prove the continuity of $t\mapsto R_n(t)$ at $t=0$ for
$n\geq 2$, we
sum once again the estimates of Lemma
\ref{lemme_perturbatif_sur_Mdbarre} and we get
$\norm{R_n(t)-R_n}_{\boL_\tau} \leq B|t|$ since the sum $\sum
m[d]D_\tau f_Z(d)$ is finite.
\end{proof}

\subsubsection{The local result}

We will apply the abstract Theorem \ref{TCLabstrait} to the
$R_n(t)$. Lemmas \ref{lemme_Rn_continu},
\ref{lemme_inversible_hors_de_1}, \ref{lemme_valeur_de_mu} and
\ref{lemme_rn_sommable} ensure that the spectral hypotheses of
this theorem are satisfied.

The only hypothesis that remains to be checked is the behavior of
the eigenvalue $\lambda(1,t)$ of $R(1,t)=\sum
R_n(t)=\hat{T_Z}(e^{itf_Z}\cdot)$. As $T_Z$ has the bounded tower
big image property and $f_Z\in L^2$, Theorem \ref{TCL_classique}
gives $\lambda(t)=1-\frac{\sigma^2}{2m(Z)} t^2+o(t^2)$.

If $\sigma^2\not=0$, we can use the abstract result Theorem
\ref{TCLabstrait} for $M(t)=t^2$ and $c=\frac{\sigma^2}{2m(Z)}$,
and we get:

\begin{prop}
\label{TCL_local} Assume $\sigma^2>0$. Then, $\forall t\in \R$, we
have
  \begin{equation*}
  \int_{T^{-n}Z\cap Z} e^{it S_nf/\sqrt{n}}\dd m \to m(Z)^2
  e^{-\sigma^2 t^2/2}.
  \end{equation*}
\end{prop}
\begin{proof}
Note that, if $t$ is fixed, $t/\sqrt{n}\in
[-\epsilon_0,\epsilon_0]$ for $n$ large enough, whence Theorem
\ref{TCLabstrait} applies. The constant $\mu$ is equal to $1/m(Z)$
according to Lemma \ref{lemme_valeur_de_mu}.

On the one hand,
  \begin{align*}
  \int T_{n,t/\sqrt{n}}1 \dd m&=\int 1_Z \hat{T}^n(e^{itS_nf/\sqrt{n}}
  1_Z)\dd m
  =\int 1_Z \circ T^n \cdot e^{itS_nf/\sqrt{n}}1_Z \dd m\\
  &=\int_{Z\cap T^{-n}Z} e^{itS_nf/\sqrt{n}}\dd m
  \end{align*}
and on the other hand,
  \begin{equation*}
  \int T_{n,t/\sqrt{n}} 1\dd m=\int
  \left(1-\frac{\sigma^2}{2\mu m(Z)}
  \frac{t\mathtwosuperior}{n}\right)^n\frac{1}{\mu}P(1)
  \dd m+o(1)
  =
  \frac{m[Z]}{\mu} e^{-\sigma^2 t^2/2} +o(1)
  \end{equation*}
using $P(1)=1_Z$, and $1/\mu=m(Z)$ whence $\mu m(Z)=1$.
\end{proof}

If $\sigma^2=0$, Theorem \ref{TCLabstrait} can not be used, but we
have nevertheless the following result:

\begin{prop}
\label{sigma_zero_donc_cocycle} Assume that $\sigma^2=0$. Then $f$
is a coboundary on $X$, i.e. it is possible to write $f=\chi\circ
T-\chi$ for some measurable function $\chi$ on $X$.
\end{prop}
\begin{proof}
Theorem \ref{TCL_classique} gives that, if $\sigma^2=0$, then $f_Z$ is a
coboundary, i.e. it can be written as $f_Z=\psi\circ T_Z-\psi$ for some
function $\psi$ on $Z$. We have to go from $f_Z$ to $f$.

We extend $\psi$ to the whole space $X$ by $\psi(x)=0$ if $x\not\in
Z$. Then, writing $\tilde{f}=f-\psi\circ T+\psi$, we have for $x\in Z$ that
  \begin{equation*}
  \sum_{0}^{\phi_Z(x)-1}\tilde{f}(T^ix)
  =f_Z(x)-\psi(T_Z x)+\psi(x)
  =0.
  \end{equation*}
We then set $\delta(x)=-\sum_{0}^{\phi_Z(x)-1} \tilde{f}(T^i x)$,
where $\phi_Z$ denotes the first return (resp. entry) time to $Z$ if
$x\in Z$ (resp. $x\not\in Z$).
Then, if $Tx\not \in Z$, we have $\tilde{f}(x)=\delta(Tx)-\delta(x)$,
while if $Tx\in Z$, then $\delta(Tx)=0$ and $\delta(x)=-\tilde{f}(x)$,
whence $\tilde{f}(x)=\delta(Tx)-\delta(x)$.
Thus, $\tilde{f}=\delta\circ T-\delta$, whence $f=(\psi+\delta)\circ
T-(\psi+\delta)$.
\end{proof}

\subsection{Proof of the Central Limit Theorem \ref{TCL_markov}}

\emph{A priori}, the constant $\sigma$ associated to $Z_q$ depends
on $q$, and should be written $\sigma_q$. We first prove that this
is not the case:
\begin{prop}
\label{constance_sigma}
$\sigma_q$ is independent of $q$.
\end{prop}
\begin{proof}
We fix $q>1$, and we want to prove that $\sigma_q=\sigma_1$. We
will write $Z=Z_q$ and $f_Z(x)=\sum_0^{\phi_Z(x)-1}f(x)$ (recall
that $f_Y$ is the same map, but on the basis $Y$ and with
the return time $\phi_Y$). We will write $\tilde{\phi}$ for the
return time from $Y$ to itself for the map $T_Z$ induced by $T$ on
$Z$. Thus, with these notations, we have for $x\in Y$ that
$f_Y(x)=\sum_{j<\tilde{\phi}(x)} f_Z(T_Z^jx)$.

Theorem \ref{TCL_classique} gives that $\sigma_1^2=\int
f_Y^2+2\int af_Y$ where $a=(I-\hat{T_Y})^{-1}(\hat{T_Y}f_Y)$ and
$\hat{T_Y}$ is the transfer operator associated to $T_Y$. In the
same way, $\sigma_q^2=\int f_Z^2+2\int bf_Z$ where
$b=(I-\hat{T_Z})^{-1}(\hat{T_Z}f_Z)$.

We define a function $c$ on $Z$ by
$c(x,i)=a(x)+\sum_{j=0}^{i-1}f_Z(T_Z^jx)$ for $x$ in the basis of
the tower and $i<\tilde{\phi}(x)$. Then, if $\psi$ is a bounded
function on $Z$, we will check that
  \begin{equation}
  \label{egalite_dans_tour}
  \int_Z c\cdot (\psi-\psi\circ T_Z) =\int_Z f_Z\cdot \psi\circ T_Z.
  \end{equation}
Indeed, note that the definition of $a$ implies that for a
function $\chi$ on $Y$, $\int_Ya\cdot (\chi-\chi \circ T_Y)=\int_Y
f_Y\cdot \chi\circ T_Y$. Thus,
  \begin{align*}
  \int_Z c\cdot (\psi-\psi\circ T_Z)
  &=\int_Y \sum_{i<\tilde{\phi}(x)} \left(a(x)+\sum_{j<i}f_Z(T_Z^jx)\right)(\psi(T_Z^i
  x)-\psi(T_Z^{i+1}x))
  \\&
  =\int_Y a(x) (\psi(x)-\psi(T_Y x)) +\sum_{j<i<\tilde{\phi}(x)}f_Z(T_Z^jx)
  (\psi(T_Z^ix)-\psi(T_Z^{i+1}x))
  \\&
  =\int_Y f_Y(x) \psi(T_Y x) +\sum_{j<\tilde{\phi}(x)-1}
  f_Z(T^j x)(\psi(T_Z^{j+1}x)-\psi(T_Y x))
  \\&
  =\int_Y \sum_{j<\tilde{\phi}(x)} f_Z(T_Z^jx) \psi(T_Y x)
  +\sum_{j<\tilde{\phi}(x)-1}
  f_Z(T^j x)(\psi(T_Z^{j+1}x)-\psi(T_Y x))
  \\&
  =\int_Y \sum_{j<\tilde{\phi}(x)-1}
  f_Z(T_Z^jx)\psi(T_Z^{j+1}x)+f_Z(T_Z^{\tilde{\phi}(x)-1}x)\psi(T_Y x)
  \\&
  =\int_Y\sum_{j<\tilde{\phi}(x)} f_Z(T_Z^j x)\psi(T_Z^{j+1}x)
  =\int_Zf_Z(x)\psi(T_Zx).
  \end{align*}

Equation \eqref{egalite_dans_tour} implies that
$(I-\hat{T_Z})c=\hat{T_Z}f_Z$. By definition of the function $b$,
we also have $(I-\hat{T_Z})b=\hat{T_Z}f_Z$, whence $c-b$ is in the
kernel of $I-\hat{T_Z}$, i.e. it is a constant $K$, and $c=b+K$.

We now compute $\sigma_1^2$, using $\int f_Z=0$:
  \begin{align*}
  \sigma_1^2
  &=\int_Yf_Y^2+2 f_Ya
  =\int_Y \left(\sum_{i<\tilde{\phi}(x)} f_Z(T_Z^ix)\right)^2 + 2\sum_{i<\tilde{\phi}(x)}
  f_Z(T_Z^ix)a(x)
  \\&
  =\int_Y \sum_{i<\tilde{\phi}(x)}f_Z(T_Z^i x)^2+2\sum_{j<i<\tilde{\phi}(x)}f_Z(T_Z^i
  x)f_Z(T_Z^j x)
  \\&\hphantom{=++}+2\sum_{i<\tilde{\phi}(x)} f_Z(T_Z^i x) \left(c(T_Z^i x)-\sum_{j<i}f_Z(T_Z^j
  x)\right)
  \\&
  =\int_Y \sum_{i<\tilde{\phi}(x)}f_Z(T_Z^i x)^2+2\sum_{i<\tilde{\phi}(x)} f_Z(T_Z^i x)
  c(T_Z^i x)
  =\int_Z f_Z(x)^2+2 f_Z(x)c(x)
  \\&
  =\int_Z f_Z(x)^2+2 f_Z(x)(b(x)+K)
  =\int_Z f_Z(x)^2+2f_Z(x)b(x)
  =\sigma_q^2.
  \end{align*}
\end{proof}

\begin{proof}[Conclusion of the proof of Theorem \ref{TCL_markov}]
Let $\sigma^2$ be given by Proposition \ref{constance_sigma}.

If $\sigma^2=0$, then Proposition \ref{sigma_zero_donc_cocycle}
ensures that $f$ can be written as $\chi\circ T-\chi$ for some
measurable function $\chi$. Hence, $S_nf/\sqrt{n}$ tends to $0$ in
probability, and the theorem is proved.

Assume now that $\sigma^2>0$. Let $t\in \R$, we want to prove that
$\int_X e^{itS_nf/\sqrt{n}} \to e^{-\sigma^2t^2/2}$. Fix
$\epsilon>0$. Let $Z=Z_q$ with $q$ large enough so that $m(Z)\geq
1-\epsilon$. Proposition \ref{TCL_local} gives $N$ such that
$\forall n\geq N$
  \begin{equation*}
  \left|\int_{T^{-n}Z\cap Z} e^{itS_nf/\sqrt{n}}-m(Z)^2e^{-\sigma^2
      t^2/2}\right| \leq \epsilon.
  \end{equation*}
Moreover,
  \begin{align*}
  \left|\int_X e^{itS_nf/\sqrt{n}}-\int_{T^{-n}Z\cap Z}
    e^{itS_nf/\sqrt{n}}\right|
  &\leq m(X-(Z\cap T^{-n}Z))\\
  &\leq m(X-Z)+m(X-T^{-n}Z)
  =2m(X-Z)
  \leq 2\epsilon.
  \end{align*}
Finally,
  \begin{equation*}
  \left|e^{-\sigma^2t^2/2}-m(Z)^2e^{-\sigma^2
      t^2/2}\right|
  \leq 1-m(Z)^2
  \leq 2\epsilon.
  \end{equation*}
The above three inequalities conclude the proof of Theorem
\ref{TCL_markov}.
\end{proof}

\subsection{Proof of Theorem \ref{thm_loi_stable_markov} on stable laws}

The proof is almost the same as that of the Central Limit Theorem. We
consider only the case $p\in (0,1)\cup (1,2)$, the case $p=2$ being
analogous using \cite{aaronson_denker:central}. For $p>1$, we can
assume without loss of generality that $f$ is centered, i.e. $\int f=0$.

We work on $Z=Z_q$, the union of the $q$ first floors of the Kakutani
tower from Section \ref{construit_Kakutani}. 
Setting $f_Z(x)=\sum_0^{\phi_Z(x)-1} f(x)$, we show as in
the proof of the central limit theorem that $\sum
\phi_Z[d]m[d]D_\tau f_Z(d)<\infty$. To apply Theorem
\ref{loi_stable_classique}, we have to estimate $m[f_Z>x]$ and
$m[f_Z<-x]$ for large $x$.

On the basis, $c_1$ is characterized by
$m[f_Y>x]=(c_1+o(1))x^{-p}L(x)$. We will show that we also have
$m[f_Z>x]=(c_1+o(1))x^{-p}L(x)$. Let $\Delta$ be the highest floor
in $Z$, and $\Lambda$ its projection on the basis $Y$. Then $f_Z$
is bounded on $Z-\Delta$ and $f_Y$ is bounded on $Y-\Lambda$,
whence it suffices to consider $\Delta$ and $\Lambda$. Let $x\in
\Lambda$ and $y=T^{q-1} x$ its image in $\Delta$. Then
$f_Z(y)=f_Y(x)+\sum_{i=0}^{q-2} f(T^i x)$. As the sum is bounded
by $C=q\norm{f_{|Z}}_\infty$, we get $f_Z(y)=f_Y(x)\pm C$. Thus,
for $x$ large enough, $m[f_Y>x-C]\geq m[f_Z>x] \geq m[f_Y>x+C]$.
As $m[f_Y>x]=(c_1+o(1))x^{-p}L(x)$, we obtain the same estimate
for $m[f_Z>x]$, using the slow variation of $L$. In the same way,
$m[f_Z<-x]=(c_2+o(1))x^{-p} L(x)$.

Using Theorem
\ref{loi_stable_classique},
we obtain that $\lambda(t)=1-a
|t|^pL(1/|t|)+ib \sgn(t)|t|^pL(1/|t|)+o(|t|^pL(1/|t|))$,
for $a=\frac{c}{m(Z)}$ and
$b=\frac{c}{m(Z)} \beta \tan\left(\frac{p\pi}{2}\right)$.

We use then the remark following Theorem \ref{TCLabstrait}, with
$M(t)=t^pL(1/t)$,
to get
  \begin{equation*}
  \norm{T_{n,t}-\frac{1}{\mu}\left( 1-c|t|^pL(1/|t|)+ic \beta
  \tan\left(\frac{p\pi}{2}\right) \sgn(t)|t|^pL(1/|t|)\right)^n P} \leq
  \epsilon(t)+\delta(n).
  \end{equation*}
If $B_n$ is defined by $B_n^p=nL(B_n)$, we have
  \begin{equation*}
  n \frac{|t|^p}{B_n^p}L\left(\frac{B_n}{|t|}\right)
  =|t|^p\frac{ L(B_n/|t|)}{L(B_n)}\to |t|^p
  \end{equation*}
since $L$ is slowly varying. Thus,
  \begin{equation*}
  T_{n,t/B_n}\to
  \frac{1}{\mu} e^{-c|t|^p +ic\beta  \tan\left(\frac{p\pi}{2}\right)
  \sgn(t)|t|^p} P
  \end{equation*}
Applying this result to the function $1$ and integrating, we get
  \begin{equation*}
  \int_{T^{-n}Z\cap Z} e^{it S_n f/B_n}\dd m
  \to m(Z)^2 e^{-c|t|^p(1-i\beta\sgn(t)
  \tan(\frac{p\pi}{2}))},
  \end{equation*}
the exponential being the characteristic function of the stable law
$X_{p,c,\beta}$.

We then let $q\to\infty$ and conclude the proof as for the Central
Limit Theorem. \qed

\bibliography{biblio}
\bibliographystyle{alpha}

\end{document}